\documentclass[12pt]{amsart}
\usepackage{amsmath,amscd,amssymb}
\usepackage[all]{xy}
\theoremstyle{plain}
\newtheorem{theorem}{Theorem}[section]
\newtheorem{proposition}[theorem]{Proposition}
\newtheorem{lemma}[theorem]{Lemma}
\newtheorem{corollary}[theorem]{Corollary}
\theoremstyle{definition}
\newtheorem{definition}[theorem]{Definition}
\newtheorem{example}[theorem]{Example}
\theoremstyle{remark}
\newtheorem{remark}[theorem]{Remark}
\newcommand{\benu}{\begin{enumerate}\renewcommand{\labelenumi}{{\rm (\roman{enumi})}}\renewcommand{\itemsep}{0pt}}
\newcommand{\eenu}{\end{enumerate}}
\newcommand{\N}{\mathbb{N}}
\newcommand{\Z}{\mathbb{Z}}

\newcommand{\T}{\mathbb{T}}

\newcommand{\cK}{{\mathcal K}}
\newcommand{\cL}{{\mathcal L}}

\newcommand{\K}{\mathbb{K}}

\newcommand{\cO}{{\mathcal O}}
\newcommand{\cT}{{\mathcal T}}
\newcommand{\cM}{{\mathcal M}}
\newcommand{\cB}{{\mathcal B}}
\newcommand{\F}{{\mathcal F}}

\newcommand{\ip}[2]{\langle{#1},{#2}\rangle}

\newcommand{\bip}[2]{\big\langle{#1},{#2}\big\rangle}

\newcommand{\motimes}{\otimes_{\rm min}}
\newcommand{\Motimes}{\otimes_{\rm max}}
\newcommand{\Ca}{$C^*$-al\-ge\-bra }
\newcommand{\Cas}{$C^*$-al\-ge\-bras }

\newcommand{\Cc}{$C^*$-cor\-re\-spon\-dence }
\newcommand{\Ccs}{$C^*$-cor\-re\-spon\-dences }
\newcommand{\Csa}{$C^*$-sub\-al\-ge\-bra }
\newcommand{\Csas}{$C^*$-sub\-al\-ge\-bras }
\DeclareMathOperator{\id}{id}
\DeclareMathOperator{\Aut}{Aut}

\DeclareMathOperator{\Ad}{Ad}

\DeclareMathOperator{\cspa}{\overline{span}}

\begin{document}
\title[On $C^*$-algebras associated with $C^*$-correspondences]
{On {\boldmath $C^*$}-algebras associated 
with {\boldmath $C^*$}-correspondences}
\author{Takeshi Katsura}
\address{
Department of Mathematical Sciences,
University of Tokyo, Komaba, Tokyo, 153-8914, JAPAN}
\curraddr{
Department of Mathematics,
University of Oregon, 
Eugene, Oregon, 97403-1222, U.S.A.}
\email{katsu@ms.u-tokyo.ac.jp}
\thanks{The author was supported in part by a Research Fellowship 
for Young Scientists of the Japan Society for the Promotion of Science.}



\begin{abstract}
We study $C^*$-al\-ge\-bras 
arising from $C^*$-cor\-re\-spon\-dences, 
which was introduced by the author. 
We prove the gauge-invariant uniqueness theorem, 
and obtain conditions for our $C^*$-al\-ge\-bras to be nuclear, exact, 
or satisfy the Universal Coefficient Theorem. 
We also obtain a 6-term exact sequence of $K$-groups 
involving the $K$-groups of our $C^*$-al\-ge\-bras. 
\end{abstract}

\maketitle

\setcounter{section}{-1}

\section{Introduction}

In \cite{Ka2}, 
we introduce a method to construct 
\Cas from $C^*$-cor\-re\-spon\-denc\-es. 
This construction is similar to 
the one of Cuntz-Pimsner algebras \cite{Pi}, 
and in fact these two constructions coincide 
when the left action of a given \Cc is injective. 
However, when the left action of a \Cc is not injective, 
our construction differs from the one in \cite{Pi}. 
Our construction of \Cas from \Ccs 
whose left actions are not injective 
is motivated by the constructions of graph algebras 
of graphs with sinks in \cite{FLR}, 
\Cas from topological graphs in \cite{Ka1}, 
and crossed products by Hilbert $C^*$-bimodules in \cite{AEE}. 
In fact, our construction generalizes all of these constructions. 
In our next paper \cite{Ka3}, 
we will explain that our \Cas have a nice property 
which crossed products by automorphisms also have. 

In this paper, 
we prove several theorems on our $C^*$-al\-ge\-bras, 
which generalize or improve known results on Cuntz-Pimsner algebras 
or other classes of $C^*$-al\-ge\-bras. 
After preliminaries of \Ccs and their representations 
in Sections \ref{SecCorres} and \ref{SecRep}, 
we give definitions of our \Cas $\cT_X$ and $\cO_X$ 
for a \Cc $X$ in Section \ref{C*algOfCor}. 
Sections \ref{SecFock} and \ref{SecCore} are preparatory sections 
for our main theorems. 
In Section \ref{SecFock}, 
we review constructions of Fock spaces and Fock representations. 
Most of the results in this section have been already known. 
In Section \ref{SecCore}, 
we analyze so-called cores. 
Main theorems can be found in Sections \ref{SecGIUT}, \ref{SecNucExact} 
and \ref{SecK}. 
In Section \ref{SecGIUT}, 
we prove the gauge-invariant uniqueness theorems of our $C^*$-al\-ge\-bras, 
which will play an important role 
in the analysis of their ideals in \cite{Ka3}. 
In Section \ref{SecNucExact}, 
we give necessary and sufficient conditions for our $C^*$-al\-ge\-bras 
to be nuclear or exact. 
In Section \ref{SecK}, 
we give a 6-term exact sequence of $K$-groups 
which seems to be helpful to compute $K$-groups of our $C^*$-al\-ge\-bras. 
We also give a sufficient condition for our \Cas 
to satisfy the Universal Coefficient Theorem of \cite{RoSc}. 

\medskip

The author would like to thank Narutaka Ozawa 
who kindly gave Example \ref{Example} 
and the proof of Proposition \ref{Fix}. 
He is also grateful 
to N. Christopher Phillips for many comments, 
and to Yasuyuki Kawahigashi for constant encouragement. 
This paper was written 
while the author was staying at University of Oregon. 
He would like to thank people there for their warm hospitality. 

\medskip

We denote by $\N=\{0,1,2,\ldots\}$ 
the set of natural numbers, 
and by $\T$ the group consisting of complex numbers 
whose absolute values are $1$. 
We use a convention that 
$\gamma(A,B)=\{\gamma(a,b)\in D\mid a\in A,b\in B\}$ 
for a map $\gamma\colon A\times B\to D$ 
such as inner products, multiplications or representations. 
We denote by $\cspa\{\cdots\}$ 
the closure of linear spans of $\{\cdots\}$. 
An ideal of a \Ca means a closed two-sided ideal.

\section{$C^*$-correspondences}\label{SecCorres}

We use \cite{Lnc2} for the general reference 
of Hilbert $C^*$-modules and $C^*$-cor\-re\-spon\-denc\-es. 

\begin{definition}\label{DefHilb}
Let $A$ be a $C^*$-algebra. 
A (right) {\em Hilbert $A$-module} $X$ is a Banach space with
a right action of the \Ca $A$ 
and an $A$-valued inner product 
$\ip{\cdot}{\cdot}_X\colon X\times X\to A$ 
satisfying certain conditions.
\end{definition}

Recall that a Hilbert $A$-module $X$ is said to be {\em full} 
if $\cspa\ip{X}{X}_X=A$. 
We do not assume that 
Hilbert $C^*$-modules $X$ are full. 
For a \Ca $A$, 
$A$ itself is a Hilbert $A$-module 
where the inner product is defined by 
$\ip{\xi}{\eta}_{A}=\xi^*\eta$, 
and the right action is multiplication. 

\begin{definition}
For Hilbert $A$-modules $X,Y$, 
we denote by $\cL(X,Y)$ 
the space of all adjointable operators from $X$ to $Y$. 
For $\xi\in X$ and $\eta\in Y$, 
the operator $\theta_{\eta,\xi}\in\cL(X,Y)$ is defined 
by $\theta_{\eta,\xi}(\zeta)=\eta\ip{\xi}{\zeta}_X\in Y$ for $\zeta\in X$. 
We define $\cK(X,Y)\subset \cL(X,Y)$ by 
$$\cK(X,Y)=\cspa\{\theta_{\eta,\xi}\in \cL(X,Y)\mid \xi\in X,\eta\in Y\}.$$
For a Hilbert $A$-module $X$, 
we set $\cL(X)=\cL(X,X)$, which is a $C^*$-algebra, and 
$\cK(X)=\cK(X,X)$, which is an ideal of $\cL(X)$. 
\end{definition}

\begin{definition}\label{DefCor}
For a \Ca $A$, 
we say that $X$ is a {\em $C^*$-cor\-re\-spon\-dence} over $A$ 
when $X$ is a Hilbert $A$-module and 
a $*$-ho\-mo\-mor\-phism $\varphi_X\colon A\to \cL(X)$ is given. 
\end{definition}

We refer to $\varphi_X$ as the left action of 
a \Cc $X$. 
A \Cc $X$ over $A$ is said to be {\em non-degenerate} if 
$\cspa(\varphi_X(A)X)=X$. 
We do not assume that \Ccs are non-degenerate. 

Let $A$ be a $C^*$-al\-ge\-bra. 
We can define a left action of the \Ca $A$ 
on the Hilbert $A$-module $A$ 
by the multiplication. 
Thus we get a \Cc over $A$, 
which is called the {\em identity correspondence} over $A$ 
and denoted by $A$. 
Note that the left action $\varphi_A$ 
of the identity correspondence $A$ 
gives an isomorphism from $A$ onto $\cK(A)\subset\cL(X)$. 

\begin{definition}\label{DefTensor}
Let $X,Y$ be \Ccs over a \Ca $A$. 
We denote by $X\odot Y$ the quotient of 
the algebraic tensor product of $X$ and $Y$ 
by the subspace generated by 
$(\xi a)\otimes \eta-\xi\otimes (\varphi_Y(a)\eta)$ 
for $\xi\in X$, $\eta\in Y$ and $a\in A$. 
We can define an $A$-valued inner product, 
right and left actions of $A$ on $X\odot Y$ by 
\begin{align*}
&\ip{\xi\otimes\eta}{\xi'\otimes\eta'}_{X\otimes Y}
=\bip{\eta}{\varphi_{Y}(\ip{\xi}{\xi'}_X)\eta'}_Y\\
&(\xi\otimes\eta)a=\xi\otimes(\eta a),\quad 
\varphi_{X\otimes Y}(a)(\xi\otimes\eta)=(\varphi_X(a)\xi)\otimes\eta,
\end{align*}
for $\xi,\xi'\in X$, $\eta,\eta'\in Y$ and $a\in A$. 
One can show that these operations are well-defined and 
extend to the completion of $X\odot Y$ 
with respect to the norm coming from the $A$-valued inner product 
defined above (see \cite[Proposition 4.5]{Lnc2}).
Thus the completion of $X\odot Y$ is a \Cc over $A$. 
This \Cc is called 
the {\em tensor product} of $X$ and $Y$, 
and denoted by $X\otimes Y$. 
\end{definition}

By definition, we have 
$$X\otimes Y=\cspa\{\xi\otimes\eta\mid \xi\in X,\eta\in Y\},$$ 
and $(\xi a)\otimes \eta=\xi\otimes (\varphi_Y(a)\eta)$ 
for $\xi\in X$, $\eta\in Y$ and $a\in A$. 

\begin{definition}
For a \Cc $X$ over a \Ca $A$ 
and $n\in\N$, 
we define a \Cc $X^{\otimes n}$ over $A$ by 
$X^{\otimes 0}=A$, $X^{\otimes 1}=X$, 
and $X^{\otimes (n+1)}=X\otimes X^{\otimes n}$ 
for $n\geq 1$. 
\end{definition}

For each $n\in\N$, 
the left action $\varphi_{X^{\otimes n}}$ of the \Cc $X^{\otimes n}$ 
will be simply denoted by $\varphi_n\colon A\to\cL(X^{\otimes n})$. 
For a positive integer $n$, we have 
$$X^{\otimes n}
=\cspa\{\xi_1\otimes\xi_2\otimes\cdots\otimes\xi_n\mid 
\xi_1,\xi_2,\ldots,\xi_n\in X\}.$$ 
Note that for positive integers $n,m$, 
there exists a natural isomorphism between 
$X^{\otimes n}\otimes X^{\otimes m}$ and $X^{\otimes (n+m)}$. 
We have such isomorphisms for $m=0$, 
but for $n=0$ we just get 
an injection $X^{\otimes 0}\otimes X^{\otimes m}\to X^{\otimes m}$. 
When $X$ is non-degenerate, 
this injection is actually an isomorphism, 
but it is not surjective in general. 

\begin{definition}
Let $n$ be a positive integer, 
and take $S\in \cL(X^{\otimes n})$. 
For each $m\in\N$, 
we define $S\otimes \id_m\in \cL(X^{\otimes (n+m)})$ 
by $(S\otimes \id_m)(\xi\otimes \eta)=S(\xi)\otimes \eta$ 
for $\xi\in X^{\otimes n}$ and $\eta\in X^{\otimes m}$. 
\end{definition}

We note that $S\otimes \id_0=S$. 
The $*$-ho\-mo\-mor\-phism 
$\cL(X^{\otimes n})\ni S\mapsto S\otimes \id_m\in \cL(X^{\otimes (n+m)})$
is injective when $\varphi_X$ is injective, 
but this is not the case in general. 
When $X$ is non-degenerate, 
we can define $S\otimes \id_n\in \cL(X^{\otimes n})$ 
for $S\in\cL(X^{\otimes 0})$ and $n\geq 1$ 
because $X^{\otimes 0}\otimes X^{\otimes n}\cong X^{\otimes n}$. 
In this case, 
we have $a\otimes \id_n=\varphi_{n}(a)$ 
for $a\in A\cong\cK(X^{\otimes 0})$. 
By abuse of notation, 
for $a\in A\cong\cK(X^{\otimes 0})$ 
we use the notation $a\otimes \id_n$ 
for denoting $\varphi_{n}(a)\in \cL(X^{\otimes n})$ 
even though $X$ is degenerate. 
Note that we cannot define $S\otimes \id_n\in \cL(X^{\otimes n})$ 
for $S\in\cL(X^{\otimes 0})$ in general. 
In other words, the $*$-ho\-mo\-mor\-phism 
$\varphi_{n}\colon A\to \cL(X^{\otimes n})$ 
need not extend to a $*$-ho\-mo\-mor\-phism $\cM(A)\to \cL(X^{\otimes n})$ 
unless $X$ is non-degenerate. 

\begin{definition}
Let us take $\xi\in X^{\otimes n}$ with $n\in\N$. 
For each $m\in\N$, 
we define an operator 
$\tau_m^n(\xi)\in\cL(X^{\otimes m},X^{\otimes (n+m)})$ by 
$$\tau_m^n(\xi)\colon X^{\otimes m}\ni\eta\mapsto 
\xi\otimes\eta\in X^{\otimes (n+m)}.$$
\end{definition}

Note that for $a\in A=X^{\otimes 0}$, 
we have $\tau_m^0(a)=\varphi_m(a)\in \cL(X^{\otimes m})$ 
for each $m\in\N$. 
Note also that $\tau_0^n\colon X^{\otimes n}\to 
\cL(X^{\otimes 0},X^{\otimes n})$ 
is an isometry onto $\cK(X^{\otimes 0},X^{\otimes n})$ 
for each $n\in\N$. 
The adjoint $\tau_m^n(\xi)^*\in\cL(X^{\otimes (n+m)},X^{\otimes m})$ 
of $\tau_m^n(\xi)$ satisfies that 
$\tau_m^n(\xi)^*(\zeta\otimes\eta)
=\varphi_m(\ip{\xi}{\zeta}_{X^{\otimes n}})\eta$ 
for $\zeta\in X^{\otimes n}$, $\eta\in X^{\otimes m}$. 
It is not difficult to see the following two lemmas. 

\begin{lemma}
For $n_1,n_2,m\in\N$ and 
$\xi_1\in X^{\otimes n_1}, \xi_2\in X^{\otimes n_2}$, 
we have 
$$\tau_{n_2+m}^{n_1}(\xi_1)\tau_{m}^{n_2}(\xi_2)
=\tau_{m}^{n_1+n_2}(\xi_1\otimes \xi_2)\qquad
\mbox{in }\cL(X^{\otimes m},X^{\otimes (n_1+n_2+m)}).$$
\end{lemma}

\begin{lemma}\label{tmn}
For $n,m\in\N$, $\xi,\eta\in X^{\otimes n}$ and $a\in A$, 
we have the following; 
\begin{align*}
\mbox{{\rm (i) }}& \tau_m^{n}(\xi)\tau_m^{n}(\eta)^*=\theta_{\xi,\eta}\otimes \id_m& 
\mbox{in }&\cL(X^{\otimes (n+m)}), \\
\mbox{{\rm (ii) }}& \tau_m^{n}(\xi)^*\tau_m^{n}(\eta)
=\varphi_m(\ip{\xi}{\eta}_{X^{\otimes n}})& 
\mbox{in }& \cL(X^{\otimes m}),\\ 
\mbox{{\rm (iii) }}& \tau_m^{n}(\xi)\varphi_m(a)=\tau_m^{n}(\xi a)& 
\mbox{in }& \cL(X^{\otimes m},X^{\otimes (n+m)}),\\ 
\mbox{{\rm (iv) }}& \varphi_{n+m}(a)\tau_m^{n}(\xi)
=\tau_m^{n}(\varphi_n(a)\xi)& 
\mbox{in }& \cL(X^{\otimes m},X^{\otimes (n+m)}). 
\phantom{Takeshi Katsura}
\end{align*}
\end{lemma}

\section{Representations of $C^*$-correspondences}\label{SecRep}

\begin{definition}\label{DefRep}
A {\em representation} of a \Cc $X$ over $A$ 
on a \Ca $B$ is 
a pair 
consisting of a $*$-ho\-mo\-mor\-phism $\pi\colon A\to B$ 
and a linear map $t\colon X\to B$ satisfying 
\benu
\item $t(\xi)^*t(\eta)=\pi\big(\ip{\xi}{\eta}_X\big)$ 
for $\xi,\eta\in X$,
\item $\pi(a)t(\xi)=t\big(\varphi_X(a)\xi\big)$ for $a\in A$, $\xi\in X$. 
\eenu
We denote by $C^*(\pi,t)$ the \Ca generated 
by the images of $\pi$ and $t$ in $B$. 
\end{definition}

A representation of a \Cc was called 
an isometric covariant representation in \cite{MS}. 
Note that for a representation $(\pi,t)$ of $X$, 
we have $t(\xi)\pi(a)=t(\xi a)$ automatically 
because the condition (i) above, 
combining with the fact that $\pi$ is a $*$-ho\-mo\-mor\-phism, implies 
$$\big\|t(\xi)\pi(a)-t(\xi a)\big\|^2
=\big\|\big(t(\xi)\pi(a)-t(\xi a)\big)^*
       \big(t(\xi)\pi(a)-t(\xi a)\big)\big\|
=0.$$
Note also that for $\xi\in X$, 
we have $\|t(\xi)\|\leq\|\xi\|_X$ 
because 
$$\|t(\xi)\|^2=\|t(\xi)^*t(\xi)\|=\|\pi(\ip{\xi}{\xi}_X)\|
\leq \|\ip{\xi}{\xi}_X\|=\|\xi\|_X^2.$$

\begin{definition}
A representation $(\pi,t)$ is said to be {\em injective}
if a $*$-ho\-mo\-mor\-phism $\pi$ is injective. 
\end{definition}

By the above computation, 
we see that $t$ is isometric 
for an injective representation $(\pi,t)$. 

\begin{definition}\label{defpsi}
For a representation $(\pi,t)$ of a \Cc $X$ on $B$, 
we define a $*$-ho\-mo\-mor\-phism $\psi_t\colon \cK(X)\to B$ 
by $\psi_t(\theta_{\xi,\eta})=t(\xi)t(\eta)^*\in B$ 
for $\xi,\eta\in X$. 
\end{definition}

For the well-definedness of a $*$-ho\-mo\-mor\-phism $\psi_t$, 
see, for example, \cite[Lemma 2.2]{KPW}.
The following lemma is easily verified. 

\begin{lemma}\label{psi}
For a representation $(\pi,t)$ of a \Cc $X$ over $A$, 
we have $\pi(a)\psi_t(k)=\psi_t(\varphi_X(a)k)$ and 
$\psi_t(k)t(\xi)=t(k \xi)$ for $a\in A$, $\xi\in X$ and $k\in\cK(X)$. 
\end{lemma}

By this lemma, we see that $\psi_t$ is injective 
for an injective representation $(\pi,t)$. 

\begin{definition}\label{Deftn}
Let $(\pi,t)$ be a representation of $X$. 
We set $t^{0}=\pi$ and $t^{1}=t$. 
For $n=2,3,\ldots$, 
we define a linear map $t^n\colon X^{\otimes n}\to C^*(\pi,t)$ 
by $t^{n}(\xi\otimes\eta)=t(\xi)t^{n-1}(\eta)$ 
for $\xi\in X$ and $\eta\in X^{\otimes (n-1)}$ . 
\end{definition}

It is routine to see that $t^n$ is well-defined and 
that $(\pi,t^n)$ is a representation of 
the \Cc $X^{\otimes n}$. 
Hence we can define 
$\psi_{t^n}\colon \cK(X^{\otimes n})\to C^*(\pi,t)$ 
by $\psi_{t^n}(\theta_{\xi,\eta})=t^n(\xi)t^n(\eta)^*$
for $\xi,\eta\in X^{\otimes n}$. 
Note that $t^n$ and $\psi_{t^n}$ are isometric 
if $(\pi,t)$ is an injective representation. 

\begin{lemma}\label{t*t}
Let $(\pi,t)$ be a representation of $X$. 
Take $\xi\in X^{\otimes n}$ and $\eta\in X^{\otimes m}$ 
for $n,m\in\N$ with $n\geq m$. 
Then we have $t^m(\eta)^*t^n(\xi)=t^{n-m}(\zeta)$ 
where $\zeta=\tau_{n-m}^m(\eta)^*\xi\in X^{\otimes (n-m)}$. 
\end{lemma}

\begin{proof}
When $m=0$, this follows from the fact that $(\pi,t^n)$ 
is a representation of the \Cc $X^{\otimes n}$. 
Let $m$ be a positive integer. 
We may assume $\xi=\eta'\otimes\zeta'$ 
for $\eta'\in X^{\otimes m}$ and $\zeta'\in X^{\otimes (n-m)}$ 
because the linear span of such elements is dense in $X^{\otimes n}$. 
We have 
\begin{align*}
t^m(\eta)^*t^n(\xi)
&=t^m(\eta)^*t^m(\eta')t^{n-m}(\zeta')\\
&=\pi(\ip{\eta}{\eta'}_{X^{\otimes m}})t^{n-m}(\zeta')\\
&=t^{n-m}(\varphi_{n-m}(\ip{\eta}{\eta'}_{X^{\otimes m}})\zeta').
\end{align*}
On the other hand, we get 
$$\tau_{n-m}^m(\eta)^*\xi=
\tau_{n-m}^m(\eta)^*(\eta'\otimes\zeta')
=\varphi_{n-m}(\ip{\eta}{\eta'}_{X^{\otimes m}})\zeta'.$$ 
We are done. 
\end{proof}

\begin{proposition}\label{cspa}
For a representation $(\pi,t)$ of $X$, 
we have 
$$C^*(\pi,t)=\cspa\{t^n(\xi)t^m(\eta)^*\mid 
\xi\in X^{\otimes n},\ \eta\in X^{\otimes m},\ n,m\in\N\}.$$
\end{proposition}

\begin{proof}
Clearly, 
the right hand side is a closed $*$-invariant linear space 
which contains the images of $\pi$ and $t$, 
and is contained in $C^*(\pi,t)$. 
Hence all we have to do is 
to check that this set is closed under the multiplication, 
and this follows from Lemma \ref{t*t}. 
\end{proof}

\section{$C^*$-algebras associated with $C^*$-correspondences}
\label{C*algOfCor}

In this section, we give definitions of the \Cas $\cT_X$ and $\cO_X$ 
for a \Cc $X$. 

\begin{definition}\label{DefTX}
Let $X$ be a \Cc over a \Ca $A$. 
We denote by $(\bar{\pi}_X,\bar{t}_X)$ 
the universal representation of $X$, 
and set $\cT_X=C^*(\bar{\pi}_X,\bar{t}_X)$. 
\end{definition}

The universal representation $(\bar{\pi}_X,\bar{t}_X)$ 
can be obtained by 
taking a direct sum of sufficiently many representations. 
By the universality, 
for every representation $(\pi,t)$ of $X$ 
we have a surjection $\rho\colon \cT_X\to C^*(\pi,t)$ 
with $\pi=\rho\circ \bar{\pi}_X$ and $t=\rho\circ \bar{t}_X$. 
This surjection will be called a natural surjection. 

\begin{definition}
For a \Cc $X$ over $A$ , 
we define an ideal $J_X$ of $A$ by 
\begin{align*}
J_X&=\varphi_X^{-1}\big(\cK(X)\big)\cap \big(\ker\varphi_X\big)^{\perp}\\
&=\{a\in A\mid \varphi_X(a)\in \cK(X)\mbox{ and } ab=0
\mbox{ for all }b\in \ker\varphi_X\}
\end{align*}
\end{definition}

Note that $J_X=\varphi_X^{-1}\big(\cK(X)\big)$ 
when $\varphi_X$ is injective. 
The ideal $J_X$ is the largest ideal 
to which the restriction of $\varphi_X$ is an injection into $\cK(X)$. 
The ideal $J_X$ has the following property. 

\begin{proposition}\label{JX}
Let $X$ be a \Cc over a \Ca $A$, 
and $(\pi,t)$ be an injective representation of $X$. 
If $a\in A$ satisfies $\pi(a)\in\psi_{t}(\cK(X))$, 
then we have $a\in J_X$ and $\pi(a)=\psi_{t}(\varphi_X(a))$. 
\end{proposition}

\begin{proof}
Take $a\in A$ with $\pi(a)\in\psi_{t}(\cK(X))$. 
Let $k\in \cK(X)$ be an element with $\pi(a)=\psi_{t}(k)$. 
For each $\xi\in X$, 
we have 
$$t(\varphi_X(a)\xi)=\pi(a)t(\xi)=\psi_t(k)t(\xi)=t(k\xi).$$
Since $t$ is injective, 
we have $\varphi_X(a)\xi=k\xi$ for every $\xi\in X$. 
This implies that $\varphi_X(a)=k\in \cK(X)$. 
Thus we get $\pi(a)=\psi_{t}(\varphi_X(a))$. 
Take $b\in \ker\varphi_X$ and we will show that $ab=0$. 
We get 
$$\pi(ab)=\pi(a)\pi(b)=\psi_{t}(\varphi_X(a))\pi(b)
=\psi_{t}(\varphi_X(a)\varphi_X(b))=0.$$
Since $\pi$ is injective, we obtain $ab=0$ as desired. 
Thus $a\in J_X$. 
\end{proof}

The above proposition motivates the following definition. 

\begin{definition}
A representation $(\pi,t)$ is said to be 
{\em covariant} 
if we have $\pi(a)=\psi_t(\varphi_X(a))$ 
for all $a\in J_X$. 
\end{definition}

\begin{definition}\label{DefOX}
For a \Cc $X$ over a \Ca $A$, 
the \Ca $\cO_X$ is defined by $\cO_X=C^*(\pi_X,t_X)$ 
where $(\pi_X,t_X)$ is the universal covariant representation of $X$. 
\end{definition}

By the universality, 
for each covariant representation $(\pi,t)$ of a \Cc $X$, 
there exists a natural surjection 
$\rho\colon \cO_X\to C^*(\pi,t)$ 
satisfying $\pi=\rho\circ\pi_X$ 
and $t=\rho\circ t_X$. 

The construction of \Cas $\cO_X$ from \Ccs $X$ 
generalizes 
both the one in \cite{Pi} for \Ccs with injective left actions 
and the one in \cite{AEE} for \Ccs coming from 
Hilbert $C^*$-bimodules. 
This is also a generalization of the construction of graph algebras 
\cite{KPRR,KPR,FLR} 
and more generally \Cas 
arising from topological graphs \cite{Ka1}. 
For the detail, see \cite{Ka2}.

\section{The Fock representation}\label{SecFock}

In this section, 
we construct a representation of a given $C^*$-cor\-re\-spon\-dence, 
which is called the Fock representation. 
The Fock representation is injective, 
and from this we get an injective covariant representation. 
Most of the results in this section 
can be found in \cite{Pi} or \cite{MS}. 
We will need them 
in Sections \ref{SecNucExact} and \ref{SecK}. 
For the convenience of the readers, 
we give complete proofs. 

\begin{definition}
The Hilbert $A$-module $\F(X)$, 
obtained as the direct sum 
of the Hilbert $A$-modules $X^{\otimes 0},X^{\otimes 1},\ldots$, 
is called the {\em Fock space} of $X$. 
\end{definition}

We consider $X^{\otimes n}$ as a submodule of $\F(X)$ 
for each $n\in\N$. 
For $n,m\in\N$, 
we consider the space $\cL(X^{\otimes n},X^{\otimes m})$ 
of adjointable operators from $X^{\otimes n}$ to $X^{\otimes m}$ 
as a subspace of $\cL(\F(X))$. 

\begin{definition}
We define a $*$-ho\-mo\-mor\-phism $\varphi_\infty\colon A\to\cL(\F(X))$ 
and a linear map $\tau_\infty\colon X\to \cL(\F(X))$ by 
$$\varphi_\infty(a)=\sum_{m=0}^\infty\varphi_m(a), \qquad
\tau_\infty(\xi)=\sum_{m=0}^\infty\tau_m^1(\xi),$$ 
for $a\in A$ and $\xi\in X$, 
where we always use the strong topology 
for the infinite sum of elements in $\cL(\F(X))$. 
\end{definition}

\begin{proposition}[{\cite[Proposition 1.3]{Pi}}]\label{Fockrep}
The pair $(\varphi_\infty,\tau_\infty)$ 
is an injective representation of $X$ on $\cL(\F(X))$. 
\end{proposition}

\begin{proof}
By taking $n=1$ in Lemma \ref{tmn} (ii) and (iv), 
we see that $(\varphi_\infty,\tau_\infty)$ 
is a representation of $X$. 
It is injective because $\varphi_0\colon A\to\cL(X^{\otimes 0})$ 
is an isomorphism onto $\cK(X^{\otimes 0})$. 
\end{proof}

This representation $(\varphi_\infty,\tau_\infty)$ is called 
the {\em Fock representation}. 
From the Fock representation $(\varphi_\infty,\tau_\infty)$, 
we can define a linear map 
$\tau_\infty^n\colon X^{\otimes n}\to \cL(\F(X))$ 
for each $n\in\N$ 
as in Definition \ref{Deftn}. 
It is easy to see that 
$\tau_\infty^n(\xi)=\sum_{m=0}^\infty\tau_m^n(\xi)$ 
for $\xi\in X^{\otimes n}$ and $n\in\N$. 

\begin{proposition}\label{pi0a}
For $a\in J_X$, 
we have 
$$\varphi_\infty(a)-\psi_{\tau_\infty}(\varphi_X(a))
=\varphi_0(a)\in \cL(X^{\otimes 0})\subset \cL(\F(X)).$$
\end{proposition}

\begin{proof}
For $\xi,\eta\in X$, 
we have 
$\psi_{\tau_\infty}(\theta_{\xi,\eta})=
\sum_{m=1}^\infty \theta_{\xi,\eta}\otimes \id_{m-1}$ 
by Lemma \ref{tmn} (i). 
Hence we have 
$\psi_{\tau_\infty}(k)=\sum_{m=1}^\infty k\otimes \id_{m-1}$ 
for all $k\in\cK(X)$. 
Therefore we obtain 
$$\varphi_\infty(a)-\psi_{\tau_\infty}(\varphi_X(a))
=\sum_{m=0}^\infty \varphi_m(a)-\sum_{m=1}^\infty \varphi_X(a)\otimes \id_{m-1}
=\varphi_0(a)$$
because $\varphi_m(a)=\varphi_X(a)\otimes \id_{m-1}$ for $m\geq 1$. 
\end{proof}

\begin{corollary}\label{cap0n=0}
If $a\in A$ satisfies $\varphi_\infty(a)\in \psi_{\tau_\infty}(\cK(X))$, 
then $a=0$. 
\end{corollary}

\begin{proof}
For $a\in A$ with 
$\varphi_\infty(a)\in \psi_{\tau_\infty}(\cK(X))$, 
we have $a\in J_X$ and 
$\varphi_\infty(a)=\psi_{\tau_\infty}(\varphi_X(a))$ 
by Proposition \ref{JX}. 
By Proposition \ref{pi0a}, 
we get $\varphi_0(a)=\varphi_\infty(a)-\psi_{\tau_\infty}(\varphi_X(a))=0$. 
Thus we obtain $a=0$ because $\varphi_0$ is injective. 
\end{proof}

The set $\F(X)J_X$ is a Hilbert $J_X$-module 
(\cite[Corollary 1.4]{Ka3}), 
and we have 
$$\cK(\F(X)J_X)=\cspa\{\theta_{\xi a,\eta}\in \cK(\F(X))
\mid \xi,\eta\in \F(X), a\in J_X \},$$ 
which is an ideal of $\cL(\F(X))$. 
We see that $k\in \cK(\F(X))$ is in $\cK(\F(X)J_X)$ 
if and only if $\ip{\xi}{k\eta}\in J_X$ 
for all $\xi,\eta\in \F(X)$ 
(see \cite[Lemma 2.6]{FMR} or \cite[Lemma 1.6]{Ka3}). 

\begin{proposition}\label{Theideal}
We have $\cK(\F(X)J_X)\subset C^*(\varphi_\infty,\tau_\infty)$. 
\end{proposition}

\begin{proof}
For $\xi\in X^{\otimes n}, \eta\in X^{\otimes m}$ 
and $a\in J_X$, 
we have 
$$\theta_{\xi a,\eta}=
\tau_\infty^n(\xi)\varphi_0(a)\tau_\infty^m(\eta)^*
=\tau_\infty^n(\xi)
 \big(\varphi_\infty(a)-\psi_{\tau_\infty}(\varphi_X(a))\big)
 \tau_\infty^m(\eta)^*
\in C^*(\varphi_\infty,\tau_\infty)$$ 
by Proposition \ref{pi0a}.
Hence $\cK(\F(X)J_X)\subset C^*(\varphi_\infty,\tau_\infty)$. 
\end{proof}

Let $\sigma\colon \cL(\F(X))\to \cL(\F(X))/\cK(\F(X)J_X)$ 
be the quotient map, 
and set $\varphi=\sigma\circ\varphi_\infty$ 
and $\tau=\sigma\circ \tau_\infty$. 
By Proposition \ref{pi0a}, 
$(\varphi,\tau)$ is a covariant representation of $X$ 
on $\cL(\F(X))/\cK(\F(X)J_X)$. 
We will see that this representation $(\varphi,\tau)$ is injective. 

\begin{lemma}\label{ktimes1}
For $n\geq 1$, 
the restriction of the $*$-ho\-mo\-mor\-phism 
$\cL(X^{\otimes n})\ni S\mapsto S\otimes \id_1\in \cL(X^{\otimes (n+1)})$
to $\cK(X^{\otimes n}J_X)$ is injective. 
\end{lemma}

\begin{proof}
Take $k\in \cK(X^{\otimes n}J_X)$ with $k\otimes \id_1=0$. 
Then for all $\xi,\xi'\in X^{\otimes n}$ 
and all $\eta,\eta'\in X$, 
we have 
$$0=
\ip{\xi\otimes \eta}{(k\otimes \id_1)(\xi'\otimes \eta')}_{X^{\otimes (n+1)}}
=\bip{\eta}{\varphi_X(\ip{\xi}{k\xi'}_{X^{\otimes n}})\eta'}_{X}.$$
Hence we have $\varphi_X(\ip{\xi}{k\xi'}_{X^{\otimes n}})=0$ 
for all $\xi,\xi'\in X^{\otimes n}$. 
Since $k\in \cK(X^{\otimes n}J_X)$, 
we have $\ip{\xi}{k\xi'}_{X^{\otimes n}}\in J_X$. 
Thus $\ip{\xi}{k\xi'}_{X^{\otimes n}}=0$ 
for all $\xi,\xi'\in X^{\otimes n}$ 
because $\varphi_X$ is injective on $J_X$. 
Therefore we get $k=0$. 
Thus the restriction of the map $S\mapsto S\otimes \id_1$ 
to $\cK(X^{\otimes n}J_X)$ is injective. 
\end{proof}

\begin{lemma}\label{Leminjrep1}
For $a\in A$, 
$\varphi_\infty(a)\in\cK(\F(X))$ 
implies $\lim_{n\to\infty}\|\varphi_n(a)\|=0$. 
\end{lemma}

\begin{proof}
For each $n\in\N$, 
let $P_n\in\cL(\F(X))$ 
be the projection onto the direct summand $X^{\otimes n}\subset \F(X)$. 
Since $\varphi_n(a)=P_n\varphi_\infty(a)P_n$, 
it suffices to show that 
$\lim_{n\to\infty}\|P_nkP_n\|=0$ 
for each $k\in\cK(\F(X))$. 
We may assume $k=\theta_{\xi,\eta}$ for $\xi,\eta\in \F(X)$ 
because the linear span of such elements is dense in $\cK(\F(X))$. 
By the same reason, 
we may assume $\xi\in X^{\otimes k}$ and 
$\eta\in X^{\otimes l}$ for some $k,l\in\N$. 
Now it is clear that we have 
$\lim_{n\to\infty}\|P_nkP_n\|=0$. 
This completes the proof. 
\end{proof}

\begin{proposition}\label{injrep}
The covariant representation $(\varphi,\tau)$ is injective. 
\end{proposition}

\begin{proof}
Take $a\in A$ with $\varphi(a)=0$. 
Then we have $\varphi_\infty(a)\in\cK(\F(X)J_X)$. 
For each $n\in\N$, 
we have 
$$\varphi_n(a)=P_n\varphi_\infty(a)P_n\in 
P_n\cK(\F(X)J_X)P_n=\cK(X^{\otimes n}J_X)$$ 
where $P_n\in\cL(\F(X))$ 
is the projection onto the direct summand $X^{\otimes n}\subset \F(X)$. 
By taking $n=0$, 
we get $a\in J_X$. 
Since $\varphi_1=\varphi_X$ is injective on $J_X$, 
we have $\|a\|=\|\varphi_1(a)\|$. 
By Lemma \ref{ktimes1}, 
we have 
$\|\varphi_n(a)\|=\|\varphi_n(a)\otimes \id_1\|=\|\varphi_{n+1}(a)\|$ 
for all positive integer $n$. 
Therefore we get $\|\varphi_n(a)\|=\|a\|$ for all $n\in\N$. 
Thus we have $a=0$ 
by Lemma \ref{Leminjrep1}. 
This proves that 
the covariant representation $(\varphi,\tau)$ is injective. 
\end{proof}

As consequences of Corollary \ref{cap0n=0} and 
Proposition \ref{injrep}, 
we have the followings. 

\begin{proposition}\label{cap=0}
The universal representation $(\bar{\pi}_X,\bar{t}_X)$ of $X$ on $\cT_X$ 
satisfies that 
$\{a\in A\mid \bar{\pi}_X(a)\in \psi_{\bar{t}_X}(\cK(X))\}=0$. 
\end{proposition}

\begin{proposition}\label{isometric}
The universal covariant representation $(\pi_X,t_X)$ of $X$ 
on $\cO_X$ is injective. 
\end{proposition}

We will see in Section \ref{SecGIUT} that 
the Fock representation $(\varphi_\infty,\tau_\infty)$ 
is the universal representation, 
and $(\varphi,\tau)$ is 
the universal covariant representation. 

Note that the \Ca $C^*(\varphi_\infty,\tau_\infty)$ 
is the augmented Cuntz-Toeplitz algebra defined in \cite{Pi}, 
and the \Ca $C^*(\varphi,\tau)$ is 
the relative Cuntz-Pimsner algebra $\cO(J_X,X)$ 
defined in \cite[Definition 2.18]{MS}.

\section{Analysis of the cores}\label{SecCore}

In this section, 
we investigate the so-called cores of \Cas $C^*(\pi,t)$ 
for representations $(\pi,t)$ of a \Cc $X$. 
Fix a \Cc $X$ over a \Ca $A$, 
and a representation $(\pi,t)$ of $X$. 

\begin{definition}\label{DefBn}
For each $n\in\N$, 
we set $B_n=\psi_{t^n}\big(\cK(X^{\otimes n})\big)\subset C^*(\pi,t)$. 
\end{definition}

Note that $B_0=\pi(A)$ and 
that $B_n\cong \cK(X^{\otimes n})$ when $(\pi,t)$ is injective. 
We can easily see the next lemma. 

\begin{lemma}\label{Bn+m}
For $n,m\in\N$ with $n\geq 1$, 
we have $\cspa(t^n(X^{\otimes n})B_m t^n(X^{\otimes n})^*)=B_{n+m}$
and $t^n(X^{\otimes n})^*B_{n+m}t^n(X^{\otimes n})\subset B_m$. 
\end{lemma}

\begin{definition}\label{DefBmn}
For $m,n\in\N$ with $m\leq n$, 
we define $B_{[m,n]}\subset C^*(\pi,t)$ by 
$B_{[m,n]}=B_m+B_{m+1}+\cdots+B_n$. 
\end{definition}

We have $B_{[n,n]}=B_n$ for each $n\in\N$. 
By the next lemma, 
we see that $B_{[m,n]}$'s are \Csas of $C^*(\pi,t)$. 

\begin{lemma}\label{kk'}
For $m,n\in\N$ with $m\leq n$, 
$k\in \cK(X^{\otimes m})$ 
and $k'\in \cK(X^{\otimes n})$, 
we have 
$\psi_{t^{m}}(k)\psi_{t^n}(k')
=\psi_{t^n}((k\otimes \id_{n-m})k')$. 
\end{lemma}

\begin{proof}
It suffices to show that 
$\psi_{t^m}(k)t^n(\xi)=t^n((k\otimes \id_{n-m})\xi)$ 
for $k\in \cK(X^{\otimes m})$ and $\xi\in X^{\otimes n}$. 
When $m=0$, 
this equation follows from the fact that $(\pi,t^n)$ 
is a representation of the \Cc $X^{\otimes n}$. 
Suppose $m\geq 1$. 
We may assume $k=\theta_{\zeta,\eta}$ for $\zeta,\eta\in X^{\otimes m}$. 
We have 
\begin{align*}
\psi_{t^m}(k)t^n(\xi)
&=t^m(\zeta)t^m(\eta)^*t^n(\xi)\\
&=t^m(\zeta)t^{n-m}\big(\tau_{n-m}^m(\eta)^*\xi\big)\\
&=t^n\big(\zeta\otimes (\tau_{n-m}^m(\eta)^*\xi)\big)\\
&=t^n\big((\tau_{n-m}^m(\zeta)\tau_{n-m}^m(\eta)^*)\xi\big)\\
&=t^n((k\otimes \id_{n-m})\xi)
\end{align*}
by Lemma \ref{t*t} and Lemma \ref{tmn} (i). 
We are done. 
\end{proof}

By the above lemma, 
$B_{[k,n]}$ is an ideal of $B_{[m,n]}$ 
for $m,k,n\in\N$ with $m\leq k\leq n$. 
In particular, 
$B_n$ is an ideal of $B_{[0,n]}$ for each $n\in\N$. 

\begin{definition}\label{DefBinf}
For $m\in\N$, 
we define a \Csa $B_{[m,\infty]}$ of $C^*(\pi,t)$ by 
$B_{[m,\infty]}=\overline{\bigcup_{n=m}^\infty B_{[m,n]}}$. 
\end{definition}

Note that the \Ca $B_{[m,\infty]}$ is 
an inductive limit of 
the increasing sequence of \Cas $\{B_{[m,n]}\}_{n=m}^\infty$. 
The \Ca $B_{[0,\infty]}$ is called 
the {\em core} of the \Ca $C^*(\pi,t)$. 
The core $B_{[0,\infty]}$ naturally arises when the \Ca $C^*(\pi,t)$ 
has an action of $\T$ called a gauge action. 

\begin{definition}
A representation $(\pi,t)$ of $X$ 
is said to {\em admit a gauge action} 
if for each $z\in\T$, 
there exists a $*$-ho\-mo\-mor\-phism 
$\beta_z\colon C^*(\pi,t)\to C^*(\pi,t)$ 
such that $\beta_z(\pi(a))=\pi(a)$ 
and $\beta_z(t(\xi))=zt(\xi)$ 
for all $a\in A$ and $\xi\in X$. 
\end{definition}

If it exists, such a $*$-ho\-mo\-mor\-phism $\beta_z$ is unique. 
By the assumptions in the definition above, 
$\beta_z$ is a $*$-automorphism for all $z\in\T$ 
and the map $\beta\colon \T\to \Aut(C^*(\pi,t))$ 
is automatically a strongly continuous homomorphism. 
By the universality, 
both the universal representation $(\bar{\pi}_X,\bar{t}_X)$ on $\cT_X$ 
and the universal covariant representation $(\pi_X,t_X)$ on $\cO_X$ 
admit gauge actions. 
We denote these actions by $\bar{\gamma}\colon \T\curvearrowright \cT_X$ 
and $\gamma\colon \T\curvearrowright \cO_X$. 
It is clear that 
for a representation $(\pi,t)$ admitting a gauge action $\beta$ 
we have $\beta_z\circ\rho=\rho\circ\bar{\gamma}_z$ 
for each $z\in\T$, 
where $\rho\colon\cT_X\to C^*(\pi,t)$ is the natural surjection. 
It is also clear that 
for a covariant representation $(\pi,t)$ 
admitting a gauge action $\beta$ 
we have $\beta_z\circ\rho=\rho\circ\gamma_z$ 
for each $z\in\T$, 
where $\rho\colon\cO_X\to C^*(\pi,t)$ is the natural surjection. 

\begin{proposition}\label{fixcspa}
When a representation $(\pi,t)$ admits a gauge action $\beta$, 
the core $B_{[0,\infty]}$ coincides with 
the fixed point algebra $C^*(\pi,t)^{\beta}$. 
\end{proposition}

\begin{proof}
Since 
$$\beta_z\big(t^n(\xi)t^m(\eta)^*\big)
=z^{n-m}t^n(\xi)t^m(\eta)^*$$
for $\xi\in X^{\otimes n}$, $\eta\in X^{\otimes m}$ 
and $z\in\T$, 
it is clear that 
$B_{[0,\infty]}\subset C^*(\pi,t)^{\beta}$. 
Take $x\in C^*(\pi,t)^{\beta}$. 
By Proposition \ref{cspa}, 
there exists a sequence $\{x_k\}_{k=0}^\infty$ 
of linear sums of elements in the form $t^n(\xi)t^m(\eta)^*$
such that $x=\lim_{k\to\infty} x_k$. 
Then we have 
$$x=\int_{\T}\beta_z(x)dz=\lim_{k\to\infty} \int_{\T}\beta_z(x_k)dz$$ 
where $dz$ is the normalized Haar measure on $\T$. 
By the above computation, we get 
$\int_{\T}\beta_z(x_k)dz\in \bigcup_{n=0}^\infty B_n$ 
for every $k$. 
Thus we have $x\in B_{[0,\infty]}$. 
We have shown that 
$B_{[0,\infty]}=C^*(\pi,t)^{\beta}$. 
\end{proof}

We are going to compute the core $B_{[0,\infty]}\subset C^*(\pi,t)$. 
To do this end, 
we need the following notation. 

\begin{definition}
For a representation $(\pi,t)$ of $X$, 
we set 
$$I'_{(\pi,t)}=\{a\in A\mid \pi(a)\in B_1=\psi_t(\cK(X))\},$$ 
which is an ideal of $A$. 
For each $n\in\N$, 
we define 
$$B'_n=\psi_{t^n}\big(\cK(X^{\otimes n}I'_{(\pi,t)})\big)
\subset B_n\subset C^*(\pi,t).$$ 
\end{definition}

\begin{proposition}\label{intersection}
For each $n\in\N$, we have $B_n\cap B_{n+1}=B'_n$. 
\end{proposition}

\begin{proof}
The case $n=0$ follows from the definition of $I'_{(\pi,t)}$. 
Let $n$ be a positive integer. 
For $a\in I'_{(\pi,t)}$ and $\xi,\eta\in X^{\otimes n}$, 
we have 
$$\psi_{t^n}(\theta_{\xi a,\eta})
=t^n(\xi a)t^n(\eta)^*
=t^n(\xi)\pi(a)t^n(\eta)^*
\in B_{n+1}$$
because $\pi(a)\in B_1$. 
Hence we get $B'_n\subset B_n\cap B_{n+1}$. 
Conversely take 
$x\in B_n\cap B_{n+1}$. 
Take $k\in \cK(X^{\otimes n})$ 
with $\psi_{t^n}(k)=x$. 
For each $\xi,\eta\in X^{\otimes n}$, 
we have 
$$\pi(\ip{\xi}{k\eta}_X)
=t^n(\xi)^*\psi_{t^n}(k)t^n(\eta)
=t^n(\xi)^*xt^n(\eta)
\in B_1$$
because $x\in B_{n+1}$. 
This implies that $\ip{\xi}{k\eta}_X\in I'_{(\pi,t)}$ 
for all $\xi,\eta\in X^{\otimes n}$. 
Hence we have $k\in \cK(X^{\otimes n}I'_{(\pi,t)})$. 
Thus we get $x=\psi_{t^n}(k)\in B'_n$. 
We have shown $B_n\cap B_{n+1}=B'_n$ for all $n\in\N$. 
\end{proof}

\begin{lemma}\label{appunit}
Let $n$ be a positive integer. 
For an approximate unit $\{u_\lambda\}$ 
of $\cK(X^{\otimes n})$ and $k\in\cK(X^{\otimes (n+1)})$, 
we have $k=\lim_{\lambda}(u_\lambda\otimes \id_1)k$. 
\end{lemma}

\begin{proof}
Clearly the equality holds for 
$k=(k'\otimes \id_1)k''\in \cK(X^{\otimes (n+1)})$ 
where $k'\in \cK(X^{\otimes n})$ and $k''\in\cK(X^{\otimes (n+1)})$. 
We will show that the linear span of such elements is 
dense in $\cK(X^{\otimes (n+1)})$. 
To do so, 
it suffices to show that the linear span of elements in the form 
$(k'\otimes \id_1)\zeta$ 
with $k'\in \cK(X^{\otimes n})$ and $\zeta\in X^{\otimes (n+1)}$ 
is dense in $X^{\otimes (n+1)}$ 
because we have 
$(k'\otimes \id_1)\theta_{\zeta,\zeta'}
=\theta_{(k'\otimes \id_1)\zeta,\zeta'}$. 
For $k'=\theta_{\xi,\xi'}$ and $\zeta=\eta\otimes \eta'$ 
with $\xi,\xi',\eta\in X^{\otimes n}$ and $\eta'\in X$, 
we have 
\begin{align*}
(k'\otimes \id_1)\zeta
&=\tau_1^n(\xi)\tau_1^n(\xi')^*(\eta\otimes \eta')\\
&=\tau_1^n(\xi)(\varphi_1(\ip{\xi'}{\eta}_{X^{\otimes n}})\eta')\\
&=\xi\otimes(\varphi_X(\ip{\xi'}{\eta}_{X^{\otimes n}})\eta')\\
&=\xi\ip{\xi'}{\eta}_{X^{\otimes n}}\otimes\eta'.
\end{align*}
Since the linear span of elements in the form 
$\xi\ip{\xi'}{\eta}_{X^{\otimes n}}$ with $\xi,\xi',\eta\in X^{\otimes n}$ 
is dense in $X^{\otimes n}$ 
and the linear span of elements in the form 
$\xi\otimes\eta'$ with $\xi\in X^{\otimes n}$ and $\eta'\in X$, 
is dense in $X^{\otimes (n+1)}$, 
we see that the linear span of elements in the form 
$(k'\otimes \id_1)\zeta$ 
with $k'\in \cK(X^{\otimes n})$ and $\zeta\in X^{\otimes (n+1)}$ 
is dense in $X^{\otimes (n+1)}$. 
We are done. 
\end{proof}

\begin{proposition}\label{Bn}
For each $n\in\N$, we have $B_{[0,n]}\cap B_{n+1}\subset B_n$. 
\end{proposition}

\begin{proof}
The assertion is obvious for $n=0$. 
We assume $n\geq 1$. 
Take $x\in B_{[0,n]}\cap B_{n+1}$. 
Choose $k\in \cK(X^{\otimes (n+1)})$ 
such that $x=\psi_{t^{n+1}}(k)$. 
For an approximate unit $\{u_\lambda\}$ 
of $\cK(X^{\otimes n})$, 
we have $k=\lim_{\lambda}(u_\lambda\otimes \id_1)k$ 
by Lemma \ref{appunit}. 
Since $\psi_{t^n}(u_\lambda)\psi_{t^{n+1}}(k)=
\psi_{t^{n+1}}\big((u_\lambda\otimes \id_1)k\big)$ 
by Lemma \ref{kk'}, 
we see 
$$x=\psi_{t^{n+1}}(k)
=\lim_{\lambda}\psi_{t^n}(u_\lambda)\psi_{t^{n+1}}(k)
=\lim_{\lambda}\psi_{t^n}(u_\lambda)x.$$
Since $B_n$ is an ideal of $B_{[0,n]}$, 
we have $x\in B_n$. 
Thus we obtain $B_{[0,n]}\cap B_{n+1}\subset B_n$. 
\end{proof}

\begin{proposition}\label{[0,n]n+1}
For each $n\in\N$, 
we have $B_{[0,n]}\cap B_{n+1}=B'_n$, 
and we get the following commutative diagram with exact rows; 
$$\begin{CD}
0 @>>> B'_n @>>> B_{[0,n]} @>>> 
B_{[0,n]}/B'_n @>>> 0\phantom{.} \\
@. @VVV @VVV  @| \\
0 @>>> B_{n+1} @>>> B_{[0,n+1]} @>>> 
B_{[0,n]}/B'_n @>>> 0. \\
\end{CD}$$ 
\end{proposition}

\begin{proof}
The former part follows from 
Proposition \ref{intersection} 
and Proposition \ref{Bn}. 
The latter part follows from the former 
and the fact $B_{[0,n+1]}=B_{[0,n]}+B_{n+1}$. 
\end{proof}

\begin{proposition}\label{0[1,n]}
For $n=1,2,\ldots,\infty$, 
we have the following short exact sequences; 
$$\begin{CD}
0 @>>> B_{[1,n]} @>>> B_{[0,n]} @>>> 
B_{0}/B'_0 @>>> 0. \\
\end{CD}$$ 
\end{proposition}

\begin{proof}
We will first prove $B_{0}\cap B_{[1,n]}=B'_0$ 
by the induction with respect to $n$. 
The case that $n=1$ follows from Proposition \ref{intersection}. 
Suppose that we have proved $B_{0}\cap B_{[1,n]}=B'_0$. 
Take $x\in B_{0}\cap B_{[1,n+1]}$. 
Choose $y\in B_{[1,n]}$ and $z\in B_{n+1}$ with $x=y+z$. 
We have $z=x-y\in B_{[0,n]}\cap B_{n+1}$. 
By Proposition \ref{Bn}, 
we have $z\in B_n$. 
Thus $x=y+z\in B_{[1,n]}$. 
Hence we have shown $B_{0}\cap B_{[1,n+1]}\subset B_{0}\cap B_{[1,n]}$. 
Since the converse inclusion is obvious, 
we get $B_{0}\cap B_{[1,n+1]}=B_{0}\cap B_{[1,n]}=B'_0$. 
Thus we obtain $B_{0}\cap B_{[1,n]}=B'_0$ 
for all positive integer $n$. 
This implies the existence of the desired short exact sequences 
for $n=1,2,\ldots$, 
because $B_{[0,n]}=B_{[1,n]}+B_0$. 
By taking inductive limits, 
we obtain the short exact sequences for $n=\infty$. 
\end{proof}

The \Csas of $\cT_X$ and $\cO_X$ 
corresponding to $B_n,B_{[m,n]}$ are denoted by 
$\bar{\cB}_n, \bar{\cB}_{[m,n]}\subset \cT_X$ 
and $\cB_n,\cB_{[m,n]}\subset \cO_X$. 
By Proposition \ref{fixcspa}
we have $\cT_X^{\,\bar{\gamma}}=\bar{\cB}_{[0,\infty]}$ and 
$\cO_X^{\,\gamma}=\cB_{[0,\infty]}$. 

\begin{proposition}\label{SplitES}
There exists a short exact sequence 
$$\begin{CD}
0@>>> \bar{\cB}_{n+1} @>>> 
\bar{\cB}_{[0,n+1]}@>>> \bar{\cB}_{[0,n]}@>>> 0, 
\end{CD}$$
which splits by the natural inclusion 
$\bar{\cB}_{[0,n]}\hookrightarrow \bar{\cB}_{[0,n+1]}$. 
\end{proposition}

\begin{proof}
This follows from Proposition \ref{[0,n]n+1} 
because Proposition \ref{cap=0} implies $I'_{(\bar{\pi}_X,\bar{t}_X)}=0$. 
\end{proof}

\begin{proposition}\label{Surj}
There exists a surjection from $\cT_{X}^{\bar{\gamma}}$ to $A$. 
\end{proposition}

\begin{proof}
This follows from Proposition \ref{0[1,n]}. 
\end{proof}

\begin{proposition}\label{ForNuc}
We get the following commutative diagram with exact rows; 
$$\begin{CD}
0 @>>> \cB_{[1,n+1]} @>>> \cB_{[0,n+1]} @>>> 
A/J_X @>>> 0\phantom{.} \\
@. @VVV @VVV  @| \\
0 @>>> \cB_{[1,\infty]} @>>> \cB_{[0,\infty]} @>>> 
A/J_X @>>> 0. \\
\end{CD}$$ 
\end{proposition}

\begin{proof}
By noting that $\cB_{0}\cong A$ and $\cB_{0}'\cong J_X$, 
this follows from Proposition \ref{0[1,n]}. 
\end{proof}

\begin{proposition}\label{ForNuc0}
We get the following commutative diagram with exact rows; 
$$\begin{CD}
0 @>>> J_X @>>> A @>>> 
A/J_X @>>> 0\phantom{.} \\
@. @VVV @VV{\pi_X}V  @| \\
0 @>>> \cB_{[1,\infty]} @>>> \cO_X^{\,\gamma} @>>> 
A/J_X @>>> 0. \\
\end{CD}$$ 
\end{proposition}

\begin{proof}
This follows from Proposition \ref{0[1,n]}. 
\end{proof}

\begin{proposition}\label{Bimod}
For a \Cc $X$ over a \Ca $A$, 
the following conditions are equivalent; 
\benu
\item the injection $\pi_X\colon A\to \cO_{X}^{\ \beta}$ is an isomorphism, 
\item we have $\cB_0\supset\cB_1$, 
\item the injection $\varphi_X\colon J_X\to \cK(X)$ is an isomorphism, 
\item the \Cc $X$ comes from a Hilbert $A$-bimodule. 
\eenu
\end{proposition}

\begin{proof}
It is clear that (i) implies (ii). 
From the condition (ii), 
we obtain $\cB_n\supset\cB_{n+1}$ for all $n\in\N$ 
by Lemma \ref{Bn+m}. 
Hence (ii) implies $\cO_{X}^{\ \beta}=\cB_0=\pi_X(A)$. 
This shows the implication (ii) $\Rightarrow$ (i). 
By setting $n=0$ in Proposition \ref{[0,n]n+1}, 
we have 
the following commutative diagram with exact rows; 
$$\begin{CD}
0 @>>> J_X @>>{\pi_X}> \cB_0 @>>> 
A/J_X @>>> 0\phantom{.} \\
@. @VV{\varphi_X}V @VVV  @| \\
0 @>>> \cK(X) @>>\psi_{t_X}> \cB_{[0,1]} @>>> 
A/J_X @>>> 0. \\
\end{CD}$$ 
From this diagram, we have the equivalence (ii)$\iff$(iii). 
Finally, the equivalence (iii)$\iff$(iv) was shown in \cite{Ka2}. 
\end{proof}

\section{The gauge-invariant uniqueness theorems}\label{SecGIUT}

In this section, 
we will give conditions 
for representations or covariant representations to be universal. 
The idea of the proof can be seen in \cite[Section 4]{Ka1} 
(and also in \cite[Section 3]{Pi}, \cite[Section 4]{FMR}). 
Let us take a \Cc $X$ over a \Ca $A$. 

\begin{proposition}\label{isomFT}
For a representation $(\pi,t)$ of $X$ 
satisfying $I'_{(\pi,t)}=0$, 
the restriction of $\rho\colon \cT_X\to C^*(\pi,t)$ 
to the fixed point algebra $\cT_X^{\,\bar{\gamma}}$ is injective. 
\end{proposition}

\begin{proof}
For $n\in\N$ 
let $B_n$ and $B_{[0,n]}$ be \Csas of $C^*(\pi,t)$ 
defined in Definition \ref{DefBn} and Definition \ref{DefBmn}. 
From the condition $I'_{(\pi,t)}=0$, 
we get the following commutative diagram with exact rows; 
$$\begin{CD}
0@>>> \bar{\cB}_{n+1} @>>> 
\bar{\cB}_{[0,n+1]}@>>> \bar{\cB}_{[0,n]}@>>> 0\\
@. @VV\rho V @VV\rho V  @VV\rho V \\
0@>>> B_{n+1} @>>> B_{[0,n+1]}@>>> B_{[0,n]}@>>> 0
\end{CD}$$
by the same argument as in Proposition \ref{SplitES}. 
Since the condition $I'_{(\pi,t)}=0$ 
implies that the representation $(\pi,t)$ is injective, 
we see that the restriction of $\rho$ to $\bar{\cB}_n$ 
is injective for all $n\in\N$. 
By using this fact and the commutative diagram above, 
we can inductively show that 
the restriction of $\rho$ to $\bar{\cB}_{[0,n]}$ is injective. 
Hence the restriction of $\rho$ to 
$\cT_{X}^{\,\bar{\gamma}}=\bar{\cB}_{[0,\infty]}$ is injective. 
\end{proof}

The following is the gauge-invariant uniqueness theorem 
for the \Ca $\cT_X$. 

\begin{theorem}\label{GIUTT}
Let $X$ be a \Cc over a \Ca $A$. 
For a representation $(\pi,t)$ 
of $X$, 
the surjection $\rho\colon \cT_X\to C^*(\pi,t)$ is an isomorphism 
if and only if $(\pi,t)$ satisfies $I'_{(\pi,t)}=0$ 
and admits a gauge action. 
\end{theorem}

\begin{proof}
We had already seen that the two conditions are necessary. 
Now suppose that a representation $(\pi,t)$ admits a gauge action $\beta$, 
and satisfies $I'_{(\pi,t)}=0$. 
Take $x\in \cT_{X}$ with $\rho(x)=0$. 
Then we have 
$$\rho\Big(\int_{\T}\bar{\gamma}_z(x^*x)dz\Big)
=\int_{\T}\rho(\bar{\gamma}_z(x^*x))dz
=\int_{\T}\beta_z(\rho(x^*x))dz
=0,$$
where $dz$ is the normalized Haar measure on $\T$. 
Since $\int_{\T}\bar{\gamma}_z(x^*x)dz\in \cT_{X}^{\,\gamma}$, 
we have $\int_{\T}\bar{\gamma}_z(x^*x)dz=0$ by Proposition \ref{isomFT}. 
This implies $x^*x=0$. 
Hence $\rho$ is injective. 
\end{proof}

\begin{proposition}\label{isomF}
For an injective covariant representation $(\pi,t)$ of $X$, 
the restriction of the surjection 
$\rho\colon \cO_X\to C^*(\pi,t)$ 
to the fixed point algebra $\cO_X^{\,\gamma}$ is injective. 
\end{proposition}

\begin{proof}
For $n\in\N$ 
let $B_n$ and $B_{[0,n]}$ be \Csas of $C^*(\pi,t)$ 
defined in Definition \ref{DefBn} and Definition \ref{DefBmn}. 
Since $\psi_{t^n}$ is injective, 
the restriction of $\rho$ 
to $\cB_n$ is an isomorphism onto $B_n$. 
It is easy to see that the restriction of $\rho$ to $\cB_{[0,n]}$ 
is a surjection onto $B_{[0,n]}$ for each $n\in\N$. 
We will show that these are injective 
by the induction with respect to $n$. 
The case that $n=0$ follows from the fact that $\pi$ is injective. 
Suppose that we had shown that the restriction of $\rho$ to $\cB_{[0,n]}$ 
is an isomorphism onto $B_{[0,n]}$. 
By Proposition \ref{JX}, 
we have $I'_{(\pi_X,t_X)}=I'_{(\pi,t)}=J_X$. 
Hence the restriction of $\rho$ 
to $\cB_n'$ is an isomorphism onto $B_n'$. 
Thus we get an isomorphism 
$\cB_{[0,n]}/\cB_n'\to B_{[0,n]}/B_n'$. 
By Proposition \ref{[0,n]n+1}
we get the following commutative diagram with exact rows; 
$$\begin{CD}
0@>>> \cB_{n+1} @>>> 
\cB_{[0,n+1]}@>>> \cB_{[0,n]}/\cB_n'@>>> 0\phantom{.}\\
@. @VV\rho V @VV\rho V  @VVV \\
0@>>> B_{n+1} @>>> B_{[0,n+1]}@>>> B_{[0,n]}/B_n'@>>> 0.
\end{CD}$$
By the 5-lemma, 
we see that the surjection $\cB_{[0,n+1]}\to B_{[0,n+1]}$ 
is an isomorphism. 
Thus we have shown that the restriction of $\rho$ 
to $\cB_{[0,n]}$ is injective for all $n\in\N$. 
Hence the restriction of $\rho$ 
to $\cO_X^{\,\gamma}=\cB_{[0,\infty]}$ is injective. 
\end{proof}

The following is the gauge-invariant uniqueness theorem 
for the \Ca $\cO_X$. 

\begin{theorem}\label{GIUT}
For a covariant representation $(\pi,t)$ of a \Cc $X$, 
the $*$-ho\-mo\-mor\-phism $\rho\colon \cO_X\to C^*(\pi,t)$ is an isomorphism 
if and only if $(\pi,t)$ is injective 
and admits a gauge action. 
\end{theorem}

\begin{proof}
The proof goes similarly as in Theorem \ref{GIUTT} 
with the help of Proposition \ref{isomF}. 
\end{proof}

When the left actions of \Ccs are injective, 
Theorem \ref{GIUT} is 
the gauge-invariant uniqueness theorem 
for Cuntz-Pimsner algebras 
which was proved in \cite[Theorem 4.1]{FMR}. 
In the case that \Ccs are defined from graphs 
with or without sinks, 
this was already proved 
in \cite[Theorem 2.1]{BHRS}. 
For \Cas arising from topological graphs, 
this was proved in \cite[Theorem 4.5]{Ka1}. 

We can apply the two gauge-invariant uniqueness theorems 
to the representations $(\varphi_\infty,\tau_\infty)$ 
and $(\varphi,\tau)$ in Section \ref{SecFock}. 

\begin{proposition}\label{GaugeActions}
Both the representation $(\varphi_\infty,\tau_\infty)$ 
and the covariant representation $(\varphi,\tau)$ 
are universal, 
that is, we have natural isomorphisms 
$C^*(\varphi_\infty,\tau_\infty)\cong \cT_X$ and 
$C^*(\varphi,\tau)\cong \cO_X$. 
\end{proposition}

\begin{proof}
To apply Theorem \ref{GIUTT} and Theorem \ref{GIUT}, 
it suffices to see that 
both of the representations $(\varphi_\infty,\tau_\infty)$ 
and $(\varphi,\tau)$ 
admit gauge actions because the other conditions 
had already been checked in Section \ref{SecFock}. 

For each $z\in\T$, define a unitary $u_z\in\cL(\F(X))$ 
by $u_z(\xi)=z^n\xi$ for $\xi\in X^{\otimes n}\subset \F(X)$ and $n\in\N$. 
It is routine to see that 
the automorphisms $\Ad u_z$ of $\cL(\F(X))$, defined by 
$\Ad u_z(x)=u_z xu_z^*$ for $x\in \cL(\F(X))$, 
give a gauge action for the representation $(\varphi_\infty,\tau_\infty)$. 
The ideal $\cK(\F(X)J_X)$ of $\cL(\F(X))$ 
is closed under the automorphisms $\Ad u_z$ for each $z\in\T$. 
Hence we can define an automorphism 
$\beta_z$ of $\cL(\F(X))/\cK(\F(X)J_X)$ 
by $\beta_z(\sigma(x))=\sigma(u_zxu_z^*)$ 
for $x\in \cL(\F(X))$ and $z\in\T$. 
It is clear that $\beta$ is a gauge action 
for the representation $(\varphi,\tau)$. 
We are done. 
\end{proof}

By Proposition \ref{GaugeActions}, 
the \Ca $\cO_X$ is isomorphic to 
the relative Cuntz-Pimsner algebras $C^*(\varphi,\tau)=\cO(J_X,X)$ 
introduced in \cite{MS} 
(cf. \cite[Theorem 2.19]{MS}). 
The isomorphism $C^*(\varphi_\infty,\tau_\infty)\cong \cT_X$ 
was already proved in \cite[Theorem 3.4]{Pi} 
under small assumption on $C^*$-cor\-re\-spon\-dences. 

The \Ca $\cO_X$ 
was defined as the largest \Ca 
among \Cas $C^*(\pi,t)$ 
generated by covariant representations $(\pi,t)$ of $X$. 
Theorem \ref{GIUT} tells us 
that we have $C^*(\pi,t)\cong\cO_X$ 
when a covariant representation $(\pi,t)$ 
satisfies two conditions; 
being injective and admitting a gauge action. 
In the next paper \cite{Ka3}, 
we will see that 
the \Ca $\cO_X$ 
can be defined as the smallest \Ca 
among \Cas $C^*(\pi,t)$ 
generated by representations $(\pi,t)$ of $X$ 
which satisfy the two conditions above; 
being injective and admitting gauge actions. 
Thus we can define $\cO_X$ 
without using the ideal $J_X$.

\section{Nuclearity and exactness}\label{SecNucExact}

In this section, 
we study when the \Cas $\cT_X$ and $\cO_X$ become nuclear or exact. 
We use the facts on nuclearity and exactness 
appeared in Appendices \ref{SecNucMap} and \ref{SecLink} 
as well as in \cite{Ws}. 

On the exactness of $\cT_X$ and $\cO_X$, 
we have the following 
which generalizes \cite[Theorem 3.1]{DS} slightly. 

\begin{theorem}[{cf. \cite[Theorem 3.1]{DS}}]\label{Exact}
For a \Cc $X$ over a \Ca $A$, 
the following conditions are equivalent;
\benu
\item $A$ is exact, 
\item $\cT_{X}^{\bar{\gamma}}$ is exact, 
\item $\cT_{X}$ is exact, 
\item $\cO_X^{\gamma}$ is exact, 
\item $\cO_{X}$ is exact. 
\eenu
\end{theorem}

\begin{proof}
Suppose that $A$ is exact. 
By Proposition \ref{cK(X)}, 
$\cK(X^{\otimes n})$ is exact for all $n\in\N$. 
By Proposition \ref{SplitES}, 
we can prove inductively that 
$\bar{\cB}_{[0,n]}\subset \cT_X^{\bar{\gamma}}$ 
is exact for all $n\in\N$ 
because exactness is closed under taking splitting extensions. 
Thus $\cT_{X}^{\bar{\gamma}}$ is exact 
because it is an inductive limit of exact $C^*$-al\-ge\-bras. 
This proves (i) $\Rightarrow$ (ii). 
The equivalences (ii)$\iff$(iii) and (iv)$\iff$(v) 
follow from Proposition \ref{Fix}. 
Since there exists a surjection $\cT_{X}\to \cO_{X}$, 
(iii) implies (v). 
Finally, (v) implies (i) because 
$\pi_X(A)\subset \cO_{X}$ is isomorphic to $A$. 
\end{proof}

On the nuclearity of $\cT_X$, 
we have the following. 

\begin{theorem}\label{TNuc}
For a \Cc $X$ over a \Ca $A$, 
the following conditions are equivalent; 
\benu
\item $A$ is nuclear, 
\item $\cT_{X}^{\bar{\gamma}}$ is nuclear, 
\item $\cT_{X}$ is nuclear. 
\eenu
\end{theorem}

\begin{proof}
In a similar way to the proof of (i) $\Rightarrow$ (ii) 
in Theorem \ref{Exact}, 
we can show that (i) implies (ii). 
The implication (ii) $\Rightarrow$ (i) follows 
from Proposition \ref{Surj}. 
Finally, Proposition \ref{Fix} gives the equivalence (ii)$\iff$(iii). 
\end{proof}

On the nuclearity of $\cO_X$, 
we have the following. 

\begin{theorem}\label{Onuc}
For a \Cc $X$ over a \Ca $A$, 
the following conditions are equivalent; 
\benu
\item $A/J_X$ is a nuclear $C^*$-al\-ge\-bra, 
and $\pi_X\colon J_X\to \cB_{[1,\infty]}$ is a nuclear map, 
\item $\pi_X\colon A\to \cO_X^{\gamma}$ is a nuclear map, 
\item $\pi_X\colon A\to \cO_X$ is a nuclear map, 
\item $\cO_X^{\gamma}$ is nuclear, 
\item $\cO_{X}$ is nuclear. 
\eenu
\end{theorem}

\begin{proof}
The equivalence (i)$\iff$(ii) is shown 
by applying Proposition \ref{NucExt}
to the diagram in Proposition \ref{ForNuc0}. 
The equivalence (ii)$\iff$(iii) follows from Proposition \ref{NucRestFix}. 
Obviously (iv) implies (ii). 
Assume (ii). 
We see that $A/J_X$ is nuclear from the equivalence (i)$\iff$(ii). 
We will prove that the embedding 
$\cB_{[0,n]}\hookrightarrow \cB_{[0,\infty]}$ 
is nuclear for all $n\in\N$ 
by the induction on $n$. 
The case $n=0$ follows from the condition (ii). 
Suppose we have shown that $\cB_{[0,n]}\hookrightarrow \cB_{[0,\infty]}$
is nuclear. 
Let us set $Y_n=\cspa(t_X(X)\cB_{[0,n]})$ 
and $Y_\infty=\cspa(t_X(X)\cB_{[0,\infty]})$. 
Then by Lemma \ref{Bn+m}, 
$Y_n$ is a Hilbert $\cB_{[0,n]}$-module 
with $\cK(Y_n)\cong \cB_{[1,n+1]}$, 
and $Y_\infty$ is a Hilbert $\cB_{[0,\infty]}$-module 
with $\cK(Y_\infty)\cong \cB_{[1,\infty]}$. 
By applying Proposition \ref{NucMapK(X)} 
to the inclusions $\cB_{[0,n]}\hookrightarrow \cB_{[0,\infty]}$ 
and $Y_n\hookrightarrow Y_\infty$, 
we see that the inclusion $\cB_{[1,n+1]}\hookrightarrow \cB_{[1,\infty]}$
is nuclear. 
Now by applying Proposition \ref{NucExt}
to the diagram in Proposition \ref{ForNuc}, 
we see that $\cB_{[0,n+1]}\hookrightarrow \cB_{[0,\infty]}$
is nuclear. 
Hence we have shown that $\cB_{[0,n]}\hookrightarrow \cB_{[0,\infty]}$
is nuclear for all $n\in\N$. 
Since $\bigcup_{n\in\N}\cB_{[0,n]}$ is dense in $\cB_{[0,\infty]}$, 
we see that the identity map $\cB_{[0,\infty]}\to \cB_{[0,\infty]}$ 
is nuclear. 
Thus $\cB_{[0,\infty]}$ is a nuclear $C^*$-al\-ge\-bra. 
This shows that (ii) implies (iv). 
Finally, the equivalence (iv)$\iff$(v) 
follows from Proposition \ref{Fix}. 
\end{proof}

We give two sufficient conditions on \Ccs $X$ 
for $\cO_{X}$ to be nuclear, 
which may be useful. 
Both of them easily follows from Theorem \ref{Onuc}. 

\begin{corollary}\label{CorOnuc1}
If $A$ is nuclear then $\cO_{X}$ is nuclear. 
\end{corollary}

\begin{corollary}\label{CorOnuc2}
If both the \Ca $A/J_X$ and the $*$-ho\-mo\-mor\-phism $\varphi_X\colon J_X\to\cK(X)$ 
are nuclear, 
then $\cO_{X}$ is nuclear. 
\end{corollary}

\begin{remark}
We can prove Corollary \ref{CorOnuc1} directly 
by showing that $\cO_{X}^{\ \gamma}$ is nuclear when $A$ is nuclear 
in a similar way to the proof of (i) $\Rightarrow$ (ii) 
in Theorem \ref{Exact}. 
\end{remark}

The converses of Corollary \ref{CorOnuc1} 
and Corollary \ref{CorOnuc2} are not true 
as the following example shows. 
We would like to thank Narutaka Ozawa 
who gave us this example. 

\begin{example}\label{Example}
Let $B$ be a nuclear $C^*$-al\-ge\-bra, 
and $D$ be a non-nuclear \Csa of $B$. 
For an integer $n$, 
we define $A_n$ by $A_n=B$ for $n>0$ and $A_n=D$ for $n\leq 0$. 
We set $A=\bigoplus_{n=-\infty}^{\infty}A_n$. 
We define an injective endomorphism $\varphi\colon A\to A$ 
so that $\varphi|_{A_0}\colon A_0\to A_1$ is a natural embedding 
and $\varphi|_{A_n}\colon A_n\to A_{n+1}$ 
is an isomorphism for a non-zero integer $n$. 
Since $D$ is not nuclear, 
the injective endomorphism $\varphi$ is not nuclear. 
Let $X$ be the \Cc over $A$ 
which is isomorphic to $A$ as Hilbert $A$-modules, 
and whose left action $\varphi_X\colon A\to\cL(X)$ 
is defined as the composition of $\varphi\colon A\to A$ 
and the isomorphism $A\cong\cK(X)\subset\cL(X)$. 
Then we have $J_X=A$ 
and the map $\varphi_X\colon J_X\to\cK(X)$ 
is not nuclear as $\varphi$ is not. 
Thus the \Cc $X$ does not satisfy the assumption of 
Corollary \ref{CorOnuc1} nor Corollary \ref{CorOnuc2}. 
However, the \Ca ${\mathcal O}_{X}$ is nuclear 
because 
the fixed point algebra ${\mathcal O}_{X}^{\ \beta}$ 
is isomorphic to the inductive limit 
$\varinjlim (A,\varphi)\cong \bigoplus_{n=-\infty}^{\infty}B$, 
which is nuclear. 
\end{example}

A Hilbert $A$-bimodule $X$ is naturally considered as 
a \Cc over $A$, 
and the \Ca $\cO_{X}$ 
is isomorphic to the crossed product $A\rtimes_X\Z$ of $A$ by $X$ 
defined in \cite[Definition 2.4]{AEE} 
(see \cite[Subsection 3.3]{Ka2}). 
We have a nice characterization of 
the nuclearity of such a $C^*$-al\-ge\-bra. 

\begin{proposition}
When a \Cc $X$ over a \Ca $A$ comes from a Hilbert $A$-bimodule, 
the \Ca $\cO_{X}$ is nuclear if and only if $A$ is nuclear. 
\end{proposition}

\begin{proof}
By Proposition \ref{Bimod}, 
we see that $\pi_X\colon A\to \cO_{X}^{\ \beta}$ is an isomorphism. 
Hence the conclusion follows from Theorem \ref{Onuc}, 
or rather Proposition \ref{Fix}. 
\end{proof}

\section{$K$-groups}\label{SecK}

The purpose of this section is 
to obtain the 6-term exact sequence of $K$-groups, 
which seems to be useful to compute the $K$-groups 
$K_0(\cO_X)$ and $K_1(\cO_X)$ of $\cO_X$. 
Mainly we follow the arguments in \cite[Section 4]{Pi}. 
There, Pimsner used $KK$-theory to obtain his 6-term exact sequence. 
For this reason, he assumed the separability of the \Cas involved. 
Here, we work directly with $K$-theory instead of using $KK$-theory, 
and obtain the 6-term exact sequence 
without the assumption of separability. 

For a \Ca $A$, we denote by $K_*(A)$ 
the $K$-group $K_0(A)\oplus K_1(A)$ of $A$ 
which has a $\Z/2\Z$-grading. 
By maps between $K$-groups, 
we mean group homomorphisms which preserve the grading. 
Thus for \Cas $A$ and $B$, 
considering maps between $K$-groups $K_*(A)\to K_*(B)$ 
is same as considering two homomorphisms $K_0(A)\to K_0(B)$ 
and $K_1(A)\to K_1(B)$. 
For a $*$-ho\-mo\-mor\-phism $\rho\colon A\to B$, 
we denote by $\rho_*$ the map $K_*(A)\to K_*(B)$ 
induced by $\rho$. 

Fix a \Cc $X$ over a \Ca $A$. 
Since we have $\cT_X\cong C^*(\varphi_\infty,\tau_\infty)$ 
by Proposition \ref{GaugeActions}, 
there exists an embedding $j\colon \cK(\F(X)J_X)\to \cT_X$ 
by Proposition \ref{Theideal}. 
Since $C^*(\varphi,\tau)\cong \cO_X$ 
by Proposition \ref{GaugeActions}, 
we have the following short exact sequence; 
$$\begin{CD}
0 @>>> \cK(\F(X)J_X)@>j>> \cT_X@>>> \cO_X@>>> 0.
\end{CD}$$

The following two propositions enable us 
to compute the $K$-groups of $\cK(\F(X)J_X)$ and $\cT_X$. 

\begin{proposition}\label{K(FJ)}
The $*$-ho\-mo\-mor\-phism $\varphi_0\colon J_X\to \cK(\F(X)J_X)$ induces 
an isomorphism $(\varphi_0)_*\colon K_*(J_X)\to K_*(\cK(\F(X)J_X))$. 
\end{proposition}

\begin{proof}
The $*$-ho\-mo\-mor\-phism $\varphi_0\colon J_X\to \cK(\F(X)J_X)$ is 
an isomorphism onto the \Csa $\cK(X^{\otimes 0}J_X)$ of $\cK(\F(X)J_X)$. 
Since $X^{\otimes 0}J_X$ is a full Hilbert $J_X$-submodule of $\F(X)J_X$, 
$\cK(X^{\otimes 0}J_X)$ is a hereditary and full \Csa of $\cK(\F(X)J_X)$. 
Hence $(\varphi_0)_*$ is an isomorphism 
by Proposition \ref{KgroupHF}. 
\end{proof}

\begin{proposition}\label{K(TX)}
The $*$-ho\-mo\-mor\-phism $\bar{\pi}_X\colon A\to \cT_X$ induces 
an isomorphism $(\bar{\pi}_X)_*\colon K_*(A)\to K_*(\cT_X)$. 
\end{proposition}

\begin{proof}
See Appendix \ref{SecKK}. 
\end{proof}

Next, we will compute $j_*\colon K_*(\cK(\F(X)J_X))\to K_*(\cT_X)$. 

\begin{definition}
We denote by $\iota\colon J_X\hookrightarrow A$ the natural embedding. 
We define a map $[X]\colon K_*(J_X)\to K_*(A)$ 
by the composition of the map 
$(\varphi_X)_*\colon K_*(J_X)\to K_*(\cK(X))$ 
induced by the restriction of $\varphi_X$ to $J_X$ and 
the map $X_*\colon K_*(\cK(X))\to K_*(A)$ 
induced by the Hilbert $A$-module $X$ 
as in Remark \ref{DefX*}. 
\end{definition}

The map $[X]\colon K_*(J_X)\to K_*(A)$ 
is same as the map induced by the element 
$(X,\varphi_X,0)$ of $KK(J_X,A)$. 
When a \Cc $X$ is defined from 
an injective $*$-ho\-mo\-mor\-phism $\varphi\colon A\to A$, 
we have $J_X=A$ and $[X]=\varphi_*$. 
For the notation in the proof of the next lemma, 
consult Appendix \ref{SecLink}.

\begin{lemma}
The composition of the two maps $[X]\colon K_*(J_X)\to K_*(A)$ 
and $(\bar{\pi}_X)_*\colon$ $K_*(A)\to K_*(\cT_X)$ 
coincides with $(\psi_{\bar{t}_X}\circ\varphi_X)_*$. 
\end{lemma}

\begin{proof}
Let $M_2(\cT_X)$ be the \Ca of two-by-two matrices 
with entries in $\cT_X$. 
For $i,j\in\{0,1\}$, 
we denote by $\iota_{ij}$ 
the natural embedding $\cT_X\to M_2(\cT_X)$ 
onto the $i,j$-component. 
By the definition of $K$-groups, 
$(\iota_{00})_*=(\iota_{11})_*$ is 
an isomorphism. 

From the maps $\bar{\pi}_X\colon A\to \cT_X$ 
and $\bar{t}_X\colon X\to \cT_X$, 
we get a $*$-ho\-mo\-mor\-phism $\rho\colon D_X\to M_2(\cT_X)$ 
such that $\rho\circ\iota_A=\iota_{11}\circ \bar{\pi}_X$ 
and $\rho\circ\iota_X=\iota_{01}\circ \bar{t}_X$. 
We have $\rho\circ\iota_{\cK(X)}=\iota_{00}\circ \psi_{\bar{t}_X}$. 
Since $X_*$ is defined as $(\iota_A)_*^{-1}\circ (\iota_{\cK(X)})_*$, 
we have 
\begin{align*}
(\bar{\pi}_X)_*\circ X_*
&=(\bar{\pi}_X)_*\circ (\iota_A)_*^{-1}\circ 
(\iota_{\cK(X)})_*\\
&=(\iota_{11})_*^{-1}\circ \rho_*\circ 
(\iota_{\cK(X)})_*\\
&=(\iota_{11})_*^{-1}\circ (\iota_{00})_*\circ (\psi_{\bar{t}_X})_*\\
&=(\psi_{\bar{t}_X})_*.
\end{align*}
Hence we get 
$$(\bar{\pi}_X)_*\circ [X]
=(\bar{\pi}_X)_*\circ X_*\circ (\varphi_X)_*
=(\psi_{\bar{t}_X})_*\circ (\varphi_X)_*
=(\psi_{\bar{t}_X}\circ \varphi_X)_*.
$$
We are done. 
\end{proof}

\begin{lemma}
The $*$-ho\-mo\-mor\-phism $\bar{\pi}_X\circ\iota\colon J_X\to\cT_X$ 
is the sum of the two $*$-ho\-mo\-mor\-phisms 
$\psi_{\bar{t}_X}\circ\varphi_X$ and $j\circ\varphi_0$. 
\end{lemma}

\begin{proof}
If we identify $\cT_X$ and $C^*(\varphi_\infty,\tau_\infty)$, 
this follows from Proposition \ref{pi0a}. 
\end{proof}

By the two lemma above, 
the map $j_*\colon K_*(\cK(\F(X)J_X))\to K_*(\cT_X)$ 
is same as the map $\iota_*-[X]\colon K_*(J_X)\to K_*(A)$ 
modulo the isomorphisms 
$(\varphi_0)_*\colon K_*(J_X)\to K_*(\cK(\F(X)J_X))$ 
and $(\bar{\pi}_X)_*\colon K_*(A)\to K_*(\cT_X)$: 
$$\begin{CD}
K_*(\cK(\F(X)J_X)) @>>{j_*}> K_*(\cT_X)\\
@AA{(\varphi_0)_*}A  @AA{(\bar{\pi}_X)_*}A\\
K_*(J_X)@>>{\iota_*-[X]}>  K_*(A).
\end{CD}$$
Thus by rewriting the 6-term exact sequence of $K$-groups 
obtained from the short exact sequence 
$$\begin{CD}
0 @>>> \cK(\F(X)J_X)@>j>> \cT_X@>>> \cO_X@>>> 0,
\end{CD}$$
we get the following. 

\begin{theorem}[{cf. \cite[Theorem 4.9]{Pi}}]\label{Kgroup}
For a \Cc $X$ over a \Ca $A$, 
we have the following exact sequence; 
$$\begin{CD}
K_0(J_X) @>>\iota_*-[X]> K_0(A) 
@>>(\pi_X)_*>  K_0(\cO_X)\phantom{.} \\
@AAA @. @VVV \\
K_1(\cO_X) @<(\pi_X)_*<< K_1(A) 
@<\iota_*-[X]<{\phantom{tk}}< K_1(J_X).
\end{CD}$$
\end{theorem}

For a \Cc $X$ over a \Ca $A$ and 
an ideal $J$ of $A$ satisfying $\varphi_X(J)\subset\cK(X)$, 
the relative Cuntz-Pimsner algebra $\cO(J,X)$ is defined 
as the quotient $C^*(\varphi_\infty,\tau_\infty)/\cK(\F(X)J)$ 
(\cite[Definition 2.18]{MS}). 
Thus we can prove the following statement 
in the same way as the proof of Theorem \ref{Kgroup}. 

\begin{proposition}
Let $X$ be a \Cc over a \Ca $A$, 
and $J$ be an ideal of $A$ with $\varphi_X(J)\subset\cK(X)$. 
Then we have the following exact sequence; 
$$\begin{CD}
K_0(J) @>>\iota_*-[X,J]> K_0(A) 
@>>\pi_*>  K_0(\cO(J,X))\phantom{,} \\
@AAA @. @VVV \\
K_1(\cO(J,X)) @<\pi_*<< K_1(A) 
@<\iota_*-[X,J]<< K_1(J),
\end{CD}$$
where $\iota\colon J\hookrightarrow A$ is the embedding, 
$\pi\colon A\to \cO(J,X)$ is the natural $*$-ho\-mo\-mor\-phism, 
and $[X,J]\colon K_*(J)\to K_*(A)$ is defined 
by $[X,J]=X_*\circ(\varphi_X|_J)_*$. 
\end{proposition}

It is not difficult to see that 
the two $*$-ho\-mo\-mor\-phisms in 
Proposition \ref{K(FJ)} and Proposition \ref{K(TX)} 
induce $KK$-equivalences between $J_X$ and $\cK(\F(X)J_X)$ 
and between $A$ and $\cT_X$ 
when the involving \Cas are separable. 
Hence by applying ``two among three principle'' 
to the short exact sequence 
$$\begin{CD}
0 @>>> \cK(\F(X)J_X)@>j>> \cT_X@>>> \cO_X@>>> 0,
\end{CD}$$
we get the following. 

\begin{proposition}\label{UCT}
Let $X$ be a separable \Cc over a separable nuclear \Ca $A$. 
If $A$ and $J_X$ satisfy 
the Universal Coefficient Theorem of \cite{RoSc}, 
then so does $\cO_X$. 
\end{proposition}

\appendix

\section{On nuclear maps}\label{SecNucMap}

In Appendices \ref{SecNucMap} and \ref{SecLink}, 
we gather the results on nuclear maps 
and linking algebras. 
We use these results in Sections \ref{SecNucExact} and \ref{SecK}. 
Most of them should be known among the specialists. 
Some results in this appendix hold 
with less assumption. 

\begin{definition}
For \Cas $A$ and $D$, 
we denote by $A\motimes D$ (resp.\ $A\Motimes D$) 
the {\em minimal (resp.\ maximal) tensor product} of $A$ and $D$, 
and by $A\ominus D$ the kernel of 
the natural surjection $\pi_{A,D}\colon A\Motimes D\to A\motimes D$. 
\end{definition}

\begin{definition}
For a $*$-ho\-mo\-mor\-phism $\varphi\colon A\to B$, 
we can define $*$-ho\-mo\-mor\-phisms 
$\varphi\motimes\id_D\colon A\motimes D\to B\motimes D$ 
and $\varphi\Motimes\id_D\colon A\Motimes D\to B\Motimes D$ such that 
$\varphi\motimes\id_D(a\otimes d)
=\varphi\Motimes\id_D(a\otimes d)
=\varphi(a)\otimes d$ 
for $a\in A$ and $d\in D$. 
Since we have the commutative diagram; 
$$\xymatrix{
A\Motimes D \ar@{>>}[d]^{\pi_{A,D}} \ar[rr]^{\varphi\Motimes\id_D} && 
B\Motimes D\phantom{,} \ar@{>>}[d]^{\pi_{B,D}} \\
A\motimes D \ar[rr]^{\varphi\motimes\id_D} && 
B\motimes D,}$$
the restriction of $\varphi\Motimes\id_D$ 
to $A\ominus D\subset A\Motimes D$ 
induces a $*$-ho\-mo\-mor\-phism 
$\varphi\ominus \id_D\colon A\ominus D\to B\ominus D$. 
\end{definition}

\begin{definition}
A $*$-ho\-mo\-mor\-phism $\varphi\colon A\to B$ 
is said to be {\em nuclear} if for all \Ca $D$, 
the $*$-ho\-mo\-mor\-phism 
$\varphi\Motimes\id_D\colon A\Motimes D\to B\Motimes D$ 
factors through the surjection 
$\pi_{A,D}\colon A\Motimes D\to A\motimes D$; 
$$\xymatrix{
A\Motimes D \ar@{>>}[d]^{\pi_{A,D}} \ar[rr]^{\varphi\Motimes\id_D} && 
B\Motimes D\phantom{.} \ar@{>>}[d]^{\pi_{B,D}} \\
A\motimes D \ar@{.>}[urr]
\ar[rr]^{\varphi\motimes\id_D} && 
B\motimes D.}$$
A \Ca $A$ is said to be {\em nuclear} 
if $\id_A\colon A\to A$ is a nuclear map. 
\end{definition}

In other words, 
a $*$-ho\-mo\-mor\-phism $\varphi\colon A\to B$ is nuclear 
if and only if $\varphi\ominus \id_D=0$ for all \Ca $D$, 
and a \Ca $A$ is nuclear 
if and only if $A\ominus D=0$ for all \Ca $D$. 

\begin{remark}
A $*$-ho\-mo\-mor\-phism is nuclear 
if and only if it has 
the completely positive approximation property 
(see \cite{Ws}). 
\end{remark}

\begin{lemma}\label{LLP}
Let 
$$\begin{CD}
0 @>>> I @>\iota>{\phantom{tk}}> A @>\pi>> B @>>> 0 
\end{CD}$$
be a short exact sequence of $C^*$-al\-ge\-bras, 
and $D$ be a $C^*$-al\-ge\-bra. 
Then the following sequence is exact; 
$$\begin{CD}
0 @>>> I\ominus D @>\iota\ominus\id_D>{\phantom{tk}}> A\ominus D 
@>\pi\ominus\id_D>> B\ominus D. 
\end{CD}$$
If there exists an injective 
nuclear $*$-ho\-mo\-mor\-phism $A\to A'$ for some \Ca $A'$, 
then $\pi\ominus \id_D$ is surjective. 
\end{lemma}

\begin{proof}
The former statement follows from the fact that 
maximal tensor products preserve short exact sequences. 
If there exists an injective 
nuclear $*$-ho\-mo\-mor\-phism $A\to A'$ for some \Ca $A'$, 
then $A$ is exact by \cite[Proposition 7.2]{Ws}. 
Since exact \Cas have Property C \cite{Ki}, 
the sequence  
$$\begin{CD}
0 @>>> I\motimes D @>\iota\motimes \id_D>{\phantom{tk}}> A\motimes D 
@>\pi\motimes \id_D>> B\motimes D @>>> 0 
\end{CD}$$
is exact (see Proposition 5.2 and Remarks 9.5.2 in \cite{Ws}). 
Hence the conclusion follows from $3\times 3$-lemma. 
\end{proof}

\begin{proposition}\label{NucExt}
Suppose that we have a following commutative diagram with exact rows; 
$$\begin{CD}
0 @>>> I @>\iota>> A @>\pi>> B @>>> 0\\ 
@. @VV\varphi_0V  @VV{\varphi}V @| @. \\ 
0 @>>> I' @>\iota'>{\phantom{tk}}> A' @>\pi'>> B @>>> 0.
\end{CD}$$
Suppose also that $\varphi$ is injective. 
Then $\varphi$ is nuclear if and only if 
both $B$ and $\varphi_0$ are nuclear. 
\end{proposition}

\begin{proof}
Take a \Ca $D$. 
By Lemma \ref{LLP}
we have the following commutative diagram with exact rows; 
$$\begin{CD}
0 @>>> I\ominus D @>\iota\ominus \id_D>> A\ominus D 
@>\pi\ominus \id_D>> B\ominus D\\ 
@. @VV\varphi_0\ominus \id_DV  @VV{\varphi\ominus \id_D}V 
@| \\ 
0 @>>> I'\ominus D @>\iota'\ominus \id_D>{\phantom{tk}}> A'\ominus D 
@>\pi'\ominus \id_D>> B\ominus D.
\end{CD}$$
Suppose that $\varphi$ is nuclear. 
By Lemma \ref{LLP}, 
the $*$-ho\-mo\-mor\-phism $\pi\ominus \id_D$ is surjective. 
Hence we have $B\ominus D=0$ for all \Ca $D$. 
We also have $\varphi_0\ominus \id_D=0$ for all \Ca $D$ 
by the diagram above. 
Thus both $B$ and $\varphi_0$ are nuclear. 
Conversely assume that both $B$ and $\varphi_0$ are nuclear. 
Then we have $\varphi\ominus \id_D=0$ for all \Ca $D$ 
by the diagram above. 
Therefore $\varphi$ is nuclear. 
We are done. 
\end{proof}

\begin{proposition}\label{NucRestrict}
Let $A$, $B$ be $C^*$-al\-ge\-bras, 
and $A_0$, $B_0$ be \Csas of $A$ and $B$, respectively. 
Let $\varphi\colon A\to B$ be a $*$-ho\-mo\-mor\-phism 
with $\varphi(A_0)\subset B_0$. 
Let $\varphi_0\colon A_0\to B_0$ be the restriction of $\varphi$. 
When $B_0$ is a hereditary \Csa of $B$, 
the nuclearity of $\varphi$ implies 
the nuclearity of $\varphi_0$. 
\end{proposition}

\begin{proof}
When $\varphi$ is nuclear, 
its restriction $\varphi'\colon A_0\to B$ is also nuclear. 
Hence for any \Ca $D$, 
the map $\varphi'\ominus\id_D\colon A_0\ominus D\to B\ominus D$ 
is $0$. 
Since $B_0$ is a hereditary \Csa of $B$, 
we see that the inclusion $\iota\colon B_0\hookrightarrow B$ 
induces an injective $*$-ho\-mo\-mor\-phism 
$\iota\Motimes\id_D\colon B_0\Motimes D\to B\Motimes D$ 
by \cite[Theorem 3.3]{Lnc1}. 
Hence the $*$-ho\-mo\-mor\-phism 
$\iota\ominus\id_D\colon B_0\ominus D\to B\ominus D$ 
is also injective. 
This shows that 
$\varphi_0\ominus\id_D\colon A_0\ominus D\to B_0\ominus D$ 
is $0$ for all \Ca $D$. 
Thus $\varphi_0$ is injective. 
\end{proof}

The following complements the proposition above. 

\begin{proposition}\label{NucRestrict2}
With the same notation in Proposition \ref{NucRestrict}, 
when $A_0$ is a hereditary and full \Csa of $A$, 
the nuclearity of $\varphi_0$ implies 
the nuclearity of $\varphi$. 
\end{proposition}

\begin{proof}
Take a \Ca $D$. 
Since $A_0$ is a hereditary and full \Csa of $A$, 
$A_0\Motimes D$ is a hereditary and full \Csa of $A\Motimes D$. 
Hence $A_0\ominus D=(A_0\Motimes D)\cap (A\ominus D)$ 
is also hereditary and full in $A\ominus D$. 
When $\varphi_0$ is nuclear, 
the $*$-ho\-mo\-mor\-phism 
$\varphi\Motimes\id_D\colon A\Motimes D\to B\Motimes D$ 
vanishes on $A_0\ominus D$. 
Thus $\varphi\Motimes\id_D$ vanishes on $A\ominus D$. 
This shows that $\varphi$ is nuclear. 
\end{proof}

The following is an immidiate consequence of 
Proposition \ref{NucRestrict} and Proposition \ref{NucRestrict2}. 

\begin{corollary}\label{FullHerNuc}
A hereditary and full \Csa $A_0$ of a \Ca $A$ is nuclear 
if and only if $A$ is nuclear. 
\end{corollary}

We also have the following. 

\begin{proposition}\label{FullHerExt}
A hereditary and full \Csa $A_0$ of a \Ca $A$ is exact 
if and only if $A$ is exact. 
\end{proposition}

\begin{proof}
Since a \Csa of an exact \Ca is exact,  
$A_0$ is exact if $A$ is exact. 
Suppose that $A_0$ is exact. 
Take a short exact sequence of $C^*$-al\-ge\-bras; 
$$\begin{CD}
0 @>>> I @>\iota>{\phantom{tk}}> B @>\pi>> D @>>> 0. 
\end{CD}$$
All we have to do is to prove $\ker(\pi\motimes\id_A)=I\motimes A$. 
Since $A_0$ is full and hereditary in $A$, 
$B\motimes A_0$ is full and hereditary in $B\motimes A$. 
Thus $\ker(\pi\motimes\id_A)$ is generated by 
its intersection with $B\motimes A_0$, 
which is $I\motimes A_0$ by the exactness of $A_0$. 
Hence we get $\ker(\pi\motimes\id_A)=I\motimes A$. 
We are done. 
\end{proof}

\begin{remark}
We can prove Proposition \ref{FullHerExt} 
by using Proposition \ref{NucRestrict2} 
together with the deep fact that a \Ca is exact if and only if 
its one (or all) faithful representation is nuclear 
due to Kirchberg \cite{Ki}. 
We can also prove Proposition \ref{FullHerExt} 
in a similar way to 
the proof of Proposition \ref{KgroupD}. 
\end{remark}

The above investigation of hereditary \Csas 
can be extended to other classes of $C^*$-sub\-al\-ge\-bras. 
In Section \ref{SecNucExact}, 
we just need the following two results. 

\begin{proposition}\label{NucRestFix}
Let $\alpha\colon G \curvearrowright A$ be an action 
of a compact group $G$ on a \Ca $A$. 
Let $\varphi\colon D\to A$ be a $*$-ho\-mo\-mor\-phism 
whose image is contained in the fixed point algebra $A^{\alpha}$ 
of $\alpha$. 
Then the restriction $\varphi_0\colon D\to A^{\alpha}$ 
is nuclear if and only if $\varphi$ is nuclear. 
\end{proposition}

\begin{proof}
Similar as the proof of Proposition \ref{NucRestrict}. 
\end{proof}

\begin{proposition}\label{Fix}
Let $\alpha\colon G \curvearrowright A$ be an action 
of a compact group $G$ on a \Ca $A$. 
Then $A$ is nuclear or exact if and only if 
the fixed point algebra $A^{\alpha}$ is also. 
\end{proposition}

\begin{proof}
For nuclearity, 
it was proved in \cite[Proposition 2]{DLRZ}. 
It was pointed out by Narutaka Ozawa that 
the technique in \cite{DLRZ} works for exactness. 
We will sketch his argument. 

When $A$ is exact, $A^{\alpha}$ is exact. 
Assume that $A^{\alpha}$ is exact. 
Take a short exact sequence of $C^*$-al\-ge\-bras; 
$$\begin{CD}
0@>>> I@>>> B@>\pi>> D@>>> 0.
\end{CD}$$
Let us take a positive element $x$ of $\ker (\pi\motimes\id_A)$. 
To derive a contradiction, 
we assume $x\notin I\motimes A$. 
Then we can find a state $\varphi$ of $B\motimes A$ 
such that $\varphi$ vanishes on $I\motimes A$ 
and $\varphi(x)>0$. 
We set $x_0=\int_G\id_{B}\motimes\alpha_z(x)dz$ 
where $dz$ is the normalized Haar measure of $G$. 
Then we see $x_0\in B\motimes A^\alpha$. 
We have 
\begin{align*}
(\pi\motimes\id_{A^\alpha})(x_0)
&=\int_G\pi\motimes\id_{A}\big(\id_{B}\motimes\alpha_z(x)\big)dz\\
&=\int_G\id_{D}\motimes\alpha_z\big(\pi\motimes\id_{A}(x)\big)dz
=0.
\end{align*}
Since $A^\alpha$ is exact, 
we have $x_0\in I\motimes A^\alpha$. 
This leads a contradiction as 
$$0=\varphi(x_0)=\int_G\varphi(\id_{B}\motimes\alpha_z(x))dz>0.$$
Therefore we have $x\in I\motimes A$ 
for all positive element $x$ of $\ker(\pi\motimes\id_A)$. 
Thus we have shown $\ker(\pi\motimes\id_A)=I\motimes A$. 
This implies that $A$ is exact. 
\end{proof}

\section{On linking algebras}\label{SecLink}

\begin{definition}
Let $A$ be a $C^*$-al\-ge\-bra 
and $X$ be a Hilbert $A$-module. 
The \Ca $\cK(X\oplus A)$ is 
called the {\em linking algebra} of $X$, 
and denoted by $D_X$. 
\end{definition}

Since $\cK(A,X)\cong X$ and $\cK(A)\cong A$ naturally, 
we have the following matrix representation of $D_X$; 
$$D_X=\bigg(\begin{array}{cc}
\cK(X)&X\\
\widetilde{X}&A
\end{array}\bigg),$$ 
where $\widetilde{X}=\cK(X,A)$ is 
the dual left Hilbert $A$-module of $X$. 
The natural embeddings are 
denoted by 
$$\iota_{\cK(X)}\colon\cK(X)\hookrightarrow D_X,\quad 
\iota_{X}\colon X\hookrightarrow D_X,\quad\mbox{and}\quad
\iota_{A}\colon A\hookrightarrow D_X.$$
Both maps $\iota_A$ and $\iota_{\cK(X)}$ 
are injective $*$-ho\-mo\-mor\-phisms onto corners of $D_X$. 
The \Csa $A$ of $D_X$ is always full, 
but $\cK(X)$ is full in $D_X$ 
only in the case that $X$ is a full Hilbert $A$-module. 

\begin{lemma}\label{LemInd}
Let $A$ be a $C^*$-al\-ge\-bra 
and $X$ be a Hilbert $A$-module. 
For separable subsets $A_0\subset A$ and $X_0\subset X$, 
there exist a separable \Csa $A_\infty\subset A$ containing $A_0$ 
and a separable closed subspace $X_\infty$ of $X$ containing $X_0$ 
such that $X_\infty$ is a Hilbert $A_\infty$-module 
by restricting the operations of $X$. 
\end{lemma}

\begin{proof}
Let $A_1$ be the \Ca generated by $A_0+\ip{X_0}{X_0}_X$. 
We set $X_1=\cspa (X_0+X_0A_0)$ which is a closed subspace of $X$. 
We inductively define families of 
separable \Csas $\{A_n\}_{n=1}^\infty$ of $A$ 
and separable closed subspaces $\{X_n\}_{n=1}^\infty$ of $X$ 
so that $A_{n+1}$ is a \Ca generated by $A_n+\ip{X_n}{X_n}_X$, 
and that $X_{n+1}=\cspa (X_n+X_nA_n)$. 
We set $A_\infty=\overline{\bigcup_{n\in\N}A_n}$ 
and $X_\infty=\overline{\bigcup_{n\in\N}X_n}$. 
Then $A_\infty$ is a separable \Csa of $A$ containing $A_0$, 
and $X_\infty$ is a separable closed subspace of $X$ containing $X_0$. 
By the construction, we have $X_\infty A_\infty\subset X_\infty$ and 
$\ip{X_\infty}{X_\infty}_X\subset A_\infty$. 
Hence $X_\infty$ is a Hilbert $A_\infty$-module. 
\end{proof}

\begin{proposition}\label{KgroupD}
For a $C^*$-al\-ge\-bra $A$ 
and a Hilbert $A$-module $X$, 
the inclusion $\iota_A\colon A\to D_X$ 
induces an isomorphism on the $K$-groups. 
\end{proposition}

\begin{proof}
When both $A$ and $X$ are separable, 
\cite[Corollary 2.6]{Br} gives us an isometry $v$ 
in the multiplier algebra $\cM(D_X\motimes\K)$ of $D_X\motimes\K$ 
such that $\varPhi\colon D_X\motimes\K\ni x\mapsto vxv^*\in A\motimes\K$ 
is an isomorphism, 
where $\K$ is the \Ca of the compact operators on the infinite-dimensional 
separable Hilbert space. 
Since the composition of the isomorphism $\varPhi$ and the inclusion 
$\iota_A\motimes \id_{\K}\colon A\motimes\K\to D_X\motimes\K$ 
induces an identity on the $K$-groups of $D_X\motimes\K$ 
(see, for example, \cite[Lemma 4.6.2]{HR}), 
the inclusion $\iota_A\motimes \id_{\K}$ 
induces an isomorphism on the $K$-groups. 
Hence the inclusion $\iota_A\colon A\to D_X$ 
also induces an isomorphism on the $K$-groups. 

Now let $A$ be a general $C^*$-al\-ge\-bra 
and $X$ be a general Hilbert $A$-module. 
By Lemma \ref{LemInd}, 
the set of the pairs $(A_\lambda,X_\lambda)$ 
consisting of separable \Csas $A_\lambda$ of $A$ 
and separable closed subspaces $X_\lambda$ of $X$ 
such that $X_\lambda$ are Hilbert $A_\lambda$-modules 
is upward directed with respect to the inclusions, 
and satisfies 
$A=\bigcup_{\lambda}A_\lambda$, $X=\bigcup_{\lambda}X_\lambda$. 
We have $A\cong\varinjlim A_\lambda$ 
and $D_X\cong\varinjlim D_{X_\lambda}$. 
By the first part of this proof, 
the inclusion $\iota_{A_\lambda}\colon A_\lambda\to D_{X_\lambda}$ 
induces an isomorphism on the $K$-groups for all $\lambda$. 
Thus the inclusion $\iota_A\colon A\to D_X$ 
also induces an isomorphism on the $K$-groups. 
\end{proof}

\begin{remark}\label{DefX*}
Let $A$ be a $C^*$-al\-ge\-bra 
and $X$ be a Hilbert $A$-module. 
By Proposition \ref{KgroupD}, 
we can define a map $X_*\colon K_*(\cK(X))\to K_*(A)$ 
by the composition of the map 
$(\iota_{\cK(X)})_*\colon K_*(\cK(X))\to K_*(D_X)$ 
and the inverse of the isomorphism 
$(\iota_{A})_*\colon K_*(A)\to K_*(D_X)$. 
This map is 
the same map as the one defined in \cite[Definition 5.1]{Ex}. 
\end{remark}

\begin{proposition}\label{KgroupHF}
Let $A$, $B$ be $C^*$-al\-ge\-bras, 
and $\iota\colon A\to B$ be an injective $*$-ho\-mo\-mor\-phism 
onto a hereditary and full \Csa of $B$. 
Then $\iota_*$ is an isomorphism 
from $K_*(A)$ to $K_*(B)$. 
\end{proposition}

\begin{proof}
The proof goes the same way 
as the proof of \cite[Corollary 2.10]{Br}
with the help of Proposition \ref{KgroupD}. 
\end{proof}

\begin{remark}
Let $A$, $B$ be strongly Morita equivalent $C^*$-al\-ge\-bras. 
Then there exists a \Ca $D$ which 
contains $A$ and $B$ as full and hereditary $C^*$-sub\-al\-ge\-bras. 
Hence we see that the $K$-groups of $A$ and $B$ are isomorphic 
by Proposition \ref{KgroupHF}, 
and that $A$ is nuclear or exact if and only if $B$ is also 
by Corollary \ref{FullHerNuc} and Proposition \ref{FullHerExt}. 
\end{remark}

We use the two propositions below 
in Section \ref{SecNucExact}. 

\begin{proposition}\label{cK(X)}
Let $A$ be a $C^*$-al\-ge\-bra 
and $X$ be a Hilbert $A$-module. 
If $A$ is nuclear or exact, 
then $\cK(X)$ is also. 
\end{proposition}

\begin{proof}
Since $A$ is a hereditary and full \Csa of $D_X$, 
if $A$ is nuclear or exact then $D_X$ is also 
by Corollary \ref{FullHerNuc} and 
Proposition \ref{FullHerExt}. 
Now the conclusion follows from the fact that 
$\cK(X)$ is a hereditary \Csa of $D_X$. 
\end{proof}

\begin{proposition}\label{NucMapK(X)}
Let $A$ and $B$ be $C^*$-al\-ge\-bras, 
$X$ be a Hilbert $A$-module, and $Y$ be a Hilbert $B$-module. 
Let $\pi\colon A\to B$ be a $*$-ho\-mo\-mor\-phism 
and $t\colon X\to Y$ be a linear map 
satisfying $\ip{t(\xi)}{t(\eta)}_Y=\pi(\ip{\xi}{\eta}_X)$ 
for $\xi,\eta\in X$. 
We can define a $*$-ho\-mo\-mor\-phism 
$\psi_t\colon \cK(X)\to \cK(Y)$ by 
$\psi_t(\theta_{\xi,\eta})=\theta_{t(\xi),t(\eta)}$ 
for $\xi,\eta\in X$. 
Then the nuclearity of $\pi$ implies the nuclearity of $\psi_t$. 
\end{proposition}

\begin{proof}
For the well-definedness of $\psi_t$, 
see \cite[Lemma 2.2]{KPW}. 
We can define a $*$-ho\-mo\-mor\-phism 
$\rho\colon D_X\to D_Y$ 
so that $\rho\circ\iota_{A}=\iota_B\circ\pi$, 
$\rho\circ\iota_X=\iota_Y\circ t$ 
and $\rho\circ\iota_{\cK(X)}=\iota_{\cK(Y)}\circ\psi_t$. 
Since $A$ is a hereditary and full \Csa of $D_X$, 
the nuclearity of $\pi$ implies the nuclearity of $\rho$ 
by Proposition \ref{NucRestrict2}. 
Since $\cK(Y)$ is a hereditary \Csa of $D_Y$, 
the nuclearity of $\rho$ implies the nuclearity of $\psi_t$ 
by Proposition \ref{NucRestrict}. 
We are done. 
\end{proof}

\section{A proof of Proposition \ref{K(TX)}}\label{SecKK}

In this appendix, we give a $K$-theoretical proof of 
Proposition \ref{K(TX)}. 
In \cite[Theorem 4.4]{Pi}, 
Pimsner used $KK$-theory to prove this proposition 
under some hypotheses, 
one of which is 
that both $A$ and $X$ are separable. 
What we will do here is 
to get rid of $KK$-theory 
from the proof of \cite[Theorem 4.4]{Pi} 
so that we can prove this proposition 
without the assumption of separability. 
We first prepare some notation and results which we will need. 

\begin{definition}
For a \Ca $A$, 
we define $SA=C_0((0,1),A)$, 
which we often consider 
as a set of functions in $C_0((-1,1),A)$ vanishing on $(-1,0]$. 
For a $*$-ho\-mo\-mor\-phism $\varphi\colon A\to B$, 
we denote by $S\varphi\colon SA\to SB$ the $*$-ho\-mo\-mor\-phism 
defined by 
$S\varphi(f)(s)=\varphi(f(s))$ for $f\in SA$ and $s\in (0,1)$. 
\end{definition}

\begin{definition}
For a \Ca $A$ and an ideal $I$ of $A$, 
we define a \Ca $D(I,A)$ by 
$$D(I,A)=\{f\in C_0((-1,1),A)\mid 
  f(s)-f(-s)\in I\mbox{ for all }s\in (-1,1)\}.$$
We denote by $\iota$ the natural embedding 
$SI\to D(I,A)$. 
\end{definition}

\begin{lemma}\label{SI=D}
The $*$-ho\-mo\-mor\-phism $\iota\colon SI\to D(I,A)$ induces 
an isomorphism $\iota_*\colon K_*(SI)\to K_*(D(I,A))$. 
\end{lemma}

\begin{proof}
Let us define a $*$-ho\-mo\-mor\-phism 
$\pi\colon D(I,A)\to C_0((-1,0],A)$ 
by the restriction. 
Then $\pi$ is surjective and 
its kernel is $SI$. 
Hence we have the following short exact sequence 
$$\begin{CD}
0@>>> SI@>{\iota}>> D(I,A)@>{\pi}>> C_0((-1,0],A)@>>> 0. 
\end{CD}$$
The conclusion follows from the 6-term exact sequence of $K$-groups 
associated with this short exact sequence together with 
the fact that $K_*(C_0((-1,0],A))=0$. 
\end{proof}

\begin{definition}\label{Defrho}
Let $A,B$ be $C^*$-al\-ge\-bras, 
and $I$ be an ideals of $A$. 
For two $*$-ho\-mo\-mor\-phisms $\rho_+,\rho_-\colon B\to A$ 
such that $\rho_+(b)-\rho_-(b)\in I$ for all $b\in B$, 
we define a $*$-ho\-mo\-mor\-phisms $\rho\colon SB\to D(I,A)$ 
by 
$$\rho(f)(s)=\begin{cases}
\rho_+(f(s))& \mbox{ if }s\geq 0\\
\rho_-(f(-s))& \mbox{ if }s\leq 0,
\end{cases}$$
for $f\in SB$. 
\end{definition}

\begin{lemma}\label{0}
When $\rho_+=\rho_-$, 
the $*$-ho\-mo\-mor\-phism $\rho\colon SB\to D(I,A)$ 
in Definition \ref{Defrho} 
induces $0$ on $K$-groups. 
\end{lemma}

\begin{proof}
When $\rho_+=\rho_-$, 
the $*$-ho\-mo\-mor\-phism 
$\rho$ factors through the $*$-ho\-mo\-mor\-phism 
$\sigma\colon C_0([0,1),A)\to D(I,A)$
defined by 
$$\sigma(f)(s)=\begin{cases}
f(s)& \mbox{ if }s\geq 0\\
f(-s)& \mbox{ if }s\leq 0,
\end{cases}$$
for $f\in C_0([0,1),A)$. 
Since $K_*(C_0([0,1),A))=0$, 
we have $\rho_*=0$. 
\end{proof}

\begin{lemma}
For $j=1,2$, 
let $A_j$ be a $C^*$-al\-ge\-bra, 
and $I_j$ be an ideal of $A_j$. 
For a $*$-ho\-mo\-mor\-phism $\varphi\colon A_1\to A_2$ 
with $\varphi(I_1)\subset I_2$, 
we can define a $*$-ho\-mo\-mor\-phism 
$D\varphi\colon D(I_1,A_1)\to D(I_2,A_2)$ 
by $D\varphi(f)(s)=\varphi(f(s))$, 
and we get a commutative diagram; 
$$\begin{CD}
SI_1 @>>{S\varphi}> SI_2\\
@VV{\iota_1}V @VV{\iota_2}V \\
D(I_1,A_1) @>>{D\varphi}> D(I_2,A_2). 
\end{CD}$$
\end{lemma}

\begin{proof}
Straightforward. 
\end{proof}

We go back to the proof of Proposition \ref{K(TX)}. 
We first treat the case that the \Cc $X$ is non-degenerate. 
Let us take a \Ca $A$ and a non-degenerate \Cc $X$. 

Let $(\varphi_{\infty},\tau_{\infty})$ 
be the Fock representation of $X$ on $\cL(\F(X))$. 
We denote by $\rho_+\colon \cT_X\to \cL(\F(X))$ 
the $*$-ho\-mo\-mor\-phism such that 
$\rho_+\circ\bar{\pi}_X=\varphi_{\infty}$ and 
$\rho_+\circ\bar{t}_X=\tau_{\infty}$.
We define a $*$-ho\-mo\-mor\-phism 
$\varphi_\infty^-\colon A\to \cL(\F(X))$ 
and a linear map $\tau_\infty^-\colon X\to \cL(\F(X))$ by 
$$\varphi_\infty^-(a)=\sum_{m=1}^\infty\varphi_m(a), \qquad
\tau_\infty^-(\xi)=\sum_{m=1}^\infty\tau_m^1(\xi).$$ 
Similarly as the proof of Proposition \ref{Fockrep}, 
we see that $(\varphi_{\infty}^-,\tau_{\infty}^-)$ 
is a representation of $X$. 
Hence there exists a $*$-ho\-mo\-mor\-phism 
$\rho_-\colon \cT_X\to \cL(\F(X))$ 
such that 
$\rho_-\circ\bar{\pi}_X=\varphi_{\infty}^-$ and 
$\rho_-\circ\bar{t}_X=\tau_{\infty}^-$.

\begin{lemma}[{\cite[Lemma 4.2]{Pi}}]\label{Lem4.2}
For every $x\in \cT_X$, 
we have $\rho_+(x)-\rho_-(x)\in\cK(\F(X))$. 
\end{lemma}

\begin{proof}
Since $\cT_X$ is generated by the image of the two maps 
$\bar{\pi}_X$ and $\bar{t}_X$, 
it suffices to show this lemma 
when $x\in \cT_X$ is in the image of these maps. 
For $a\in A$, 
we have 
$$\rho_+(\bar{\pi}_X(a))-\rho_-(\bar{\pi}_X(a))=\varphi_0(a)\in \cK(\F(X)),$$ 
and for $\xi\in X$, 
we have 
$$\rho_+(\bar{t}_X(\xi))-\rho_-(\bar{t}_X(\xi))=\tau_0^1(\xi)\in \cK(\F(X)).$$ 
We are done. 
\end{proof}

Let us set $D=D\big(\cK(\F(X)),\cL(\F(X))\big)$. 
By Lemma \ref{Lem4.2}, we can define 
a $*$-ho\-mo\-mor\-phism $\rho\colon S\cT_X\to D$ by 
$$\rho(f)(s)=\begin{cases}
\rho_+(f(s))& \mbox{ if }s\geq 0\\
\rho_-(f(-s))& \mbox{ if }s\leq 0.
\end{cases}$$

\begin{lemma}\label{Spi0}
The $*$-ho\-mo\-mor\-phism $S\varphi_0\colon SA\to D$ 
induces an isomorphism on the $K$-groups. 
\end{lemma}

\begin{proof}
This follows from the fact that 
$\varphi_0\colon A\to \cK(\F(X))$ 
is an injection onto a hereditary and full \Csa of $\cK(\F(X))$ 
with the help of 
Proposition \ref{KgroupHF} and Lemma \ref{SI=D}. 
\end{proof}

\begin{proposition}\label{Inverse1}
The composition of $S\bar{\pi}_X\colon SA\to S\cT_X$ 
and $\rho\colon S\cT_X\to D$ induces an isomorphism on the $K$-groups. 
\end{proposition}

\begin{proof}
Since we have 
$\rho_+\circ \bar{\pi}_X=\varphi_0+\rho_-\circ \bar{\pi}_X$, 
we can see that the composition $\rho\circ S\bar{\pi}_X$ 
induces the same map as $S\varphi_0$ 
with the help of Lemma \ref{0}. 
Hence the proof completes by Lemma \ref{Spi0}. 
\end{proof}

Proposition \ref{Inverse1} implies 
that $\rho_*$ is ``the left inverse'' of the map 
$(S\bar{\pi}_X)_*\colon K_*(SA)\to K_*(S\cT_X)$ 
modulo the isomorphism $(S\varphi_0)_*$. 
We will show that it is also ``the right inverse''. 
To this end, we first ``shift'' 
the $*$-ho\-mo\-mor\-phism $S\bar{\pi}_X\colon SA\to S\cT_X$ 
along the $*$-ho\-mo\-mor\-phism $S\varphi_0\colon SA\to D$ 
(see Lemma \ref{shift}). 

\begin{definition}
For each $n\in\N$, we set 
$Y_n=\cspa(\bar{t}_X^n(X^{\otimes n})\cT_X)\subset \cT_X$, 
which is naturally a Hilbert $\cT_X$-module. 
We denote by $Y$ the direct sum of 
the Hilbert $\cT_X$-modules $\{Y_n\}_{n=0}^\infty$. 
\end{definition}

\begin{remark}
The Hilbert $\cT_X$-module $Y$ 
is isomorphic to the interior tensor product 
of the Hilbert $A$-module $\F(X)$ 
and the Hilbert $\cT_X$-module $\cT_X$ 
with the $*$-ho\-mo\-mor\-phism $\bar{\pi}_X\colon A\to \cT_X$. 
\end{remark}

The linear maps $\bar{t}_X^n\colon X^{\otimes n}\to Y_n$ 
extend a linear map $\bar{t}_X^\bullet\colon \F(X)\to Y$. 
By the definition, 
we get $Y=\cspa (\bar{t}_X^\bullet(\F(X))\cT_X)$. 
We also have 
$\ip{\bar{t}_X^\bullet(\xi)}{\bar{t}_X^\bullet(\eta)}_Y
=\bar{\pi}_X\big(\ip{\xi}{\eta}_{\F(X)}\big)$ 
for all $\xi,\eta\in \F(X)$. 

\begin{definition}
We define a $*$-ho\-mo\-mor\-phism 
$\varPhi\colon\cL(\F(X))\ni T\mapsto \varPhi(T)\in \cL(Y)$ by 
$$\varPhi(T)\big(\bar{t}_X^\bullet(\xi)x\big)=\bar{t}_X^\bullet(T(\xi))x
\quad\mbox{ for }\xi\in \F(X)\mbox{ and }x\in \cT_X.$$ 
\end{definition}

It is not difficult to see that $\varPhi$ is well-defined. 

\begin{lemma}\label{JustComp}
We have $\varPhi\big(\cK(\F(X))\big)\subset \cK(Y)$. 
\end{lemma}

\begin{proof}
This follows from the fact that 
$\varPhi(\theta_{\xi,\eta})
=\theta_{\bar{t}_X^\bullet(\xi),\bar{t}_X^\bullet(\eta)}$ 
for $\xi,\eta\in \F(X)$, 
which is easily verified. 
\end{proof}

We define $\widetilde{D}=D\big(\cK(Y),\cL(Y)\big)$. 
By Lemma \ref{JustComp}, 
we can define a $*$-ho\-mo\-mor\-phism $D\varPhi\colon D\to \widetilde{D}$. 
Since we assume that $X$ is non-degenerate, 
we have $Y_0=\cT_X$. 
Hence the natural isomorphism $\cT_X\cong\cK(Y_0)\subset\cK(Y)$ 
gives us 
a $*$-ho\-mo\-mor\-phism $\widetilde{\varphi}_0\colon \cT_X\to\cK(Y)$. 

\begin{lemma}\label{Stpi0}
The $*$-ho\-mo\-mor\-phism 
$S\widetilde{\varphi}_0\colon S\cT_X\to \widetilde{D}$ 
induces an isomorphism on the $K$-groups. 
\end{lemma}

\begin{proof}
Similar as the proof of Lemma \ref{Spi0}. 
\end{proof}

\begin{lemma}\label{shift}
We have the following commutative diagram; 
$$\begin{CD}
SA @>{S\bar{\pi}_X}>> S\cT_X \\
@VV{S \varphi_0}V @VV{S\widetilde{\varphi}_0}V \\
D @>{D\varPhi}>> \widetilde{D}, 
\end{CD}$$
\end{lemma}

\begin{proof}
Straightforward.
\end{proof}

\begin{proposition}\label{Inverse2}
The composition of $\rho\colon S\cT_X\to D$ 
and $D\varPhi\colon D\to \widetilde{D}$ 
induces an isomorphism on the $K$-groups. 
\end{proposition}

\begin{proof}
We set $\pi=\varPhi\circ\varphi_\infty\colon A\to\cL(Y)$. 
For each $s\in [0,1]$, 
we define a linear map $t_s\colon X\to \cL(Y)$ by 
$$t_s(\xi)=s\widetilde{\varphi}_0(\bar{t}_X(\xi))
 +\sqrt{1-s^2}\varPhi(\tau_0^1(\xi))
 +\varPhi(\tau_{\infty}^-(\xi))$$
It is routine to check that the pair $(\pi,t_s)$ is a representation of $X$. 
Thus we get a $*$-ho\-mo\-mor\-phism $\rho_s\colon \cT_X\to \cL(Y)$ 
such that $\rho_s\circ\bar{\pi}_X=\pi$ and 
$\rho_s\circ\bar{t}_X=t_s$ 
for each $s$. 
We have $\rho_0=\varPhi\circ\rho_+$ 
because $t_0=\varPhi\circ\tau_{\infty}$. 
We also have $\rho_1=\widetilde{\varphi}_0+\varPhi\circ\rho_-$ 
because 
$t_1=\widetilde{\varphi}_0\circ\bar{t}_X+\varPhi\circ\tau_{\infty}^-$ 
and $\pi=\widetilde{\varphi}_0\circ \bar{\pi}_X+\varPhi\circ\varphi_\infty^-$. 
For $\xi\in X$ and $s\in [0,1]$, 
we have $t_s(\xi)-\varPhi(\tau_{\infty}^-(\xi))\in\cK(Y)$ 
because $\widetilde{\varphi}_0(\bar{t}_X(\xi)), 
\varPhi(\tau_0^1(\xi))\in \cK(Y)$. 
Since we have 
$\pi(a)-\varPhi(\varphi_{\infty}^-(a))
=\widetilde{\varphi}_0(\bar{\pi}_X(a))\in \cK(Y)$, 
we can prove $\rho_s(x)-\varPhi(\rho_-(x))\in\cK(Y)$ 
for all $x\in \cT_X$ and $s\in [0,1]$ 
in a similar way to the proof of Lemma \ref{Lem4.2}. 
Hence we can see that the composition 
of $D\varPhi\circ\rho$ 
is homotopic to the $*$-ho\-mo\-mor\-phism $S\cT_X\to \widetilde{D}$ 
defined from the two $*$-ho\-mo\-mor\-phisms 
$S\widetilde{\varphi}_0+S\varPhi\circ\rho_-$ 
and $S\varPhi\circ\rho_-$. 
By Lemma \ref{0}, 
we see that $D\varPhi\circ\rho$ 
induces the same map as $S\widetilde{\varphi}_0$. 
Hence the proof completes by Lemma \ref{Stpi0}. 
\end{proof}

Combining all the results above, 
we obtain that the composition of 
the map $\rho_*\colon K_*(S\cT_X)\to K_*(D)$ 
and the isomorphism 
$(S\varphi_0)_*^{-1}\colon K_*(D)\to K_*(SA)$ 
gives the inverse of the map 
$(S\bar{\pi}_X)_*\colon K_*(SA)\to K_*(S\cT_X)$. 
Hence we have shown that 
$(\bar{\pi}_X)_*\colon K_*(A)\to K_*(\cT_X)$ is an isomorphism 
when the \Cc $X$ is non-degenerate. 
We will see that this is the case 
for general $C^*$-cor\-re\-spon\-dences. 

Let us take a \Cc $X$ over a \Ca $A$. 
We define 
$$T=\cspa(\bar{\pi}_X(A)\cT_X\bar{\pi}_X(A))$$ 
which is the hereditary \Csa of $\cT_X$ generated by $\bar{\pi}_X(A)$. 
Since the ideal generated by $\bar{\pi}_X(A)$ is $\cT_X$, 
Proposition \ref{KgroupHF} shows that 
the inclusion $T\hookrightarrow \cT_X$ 
induces an isomorphism on the $K$-groups. 
Hence to prove that the $*$-ho\-mo\-mor\-phism 
$\bar{\pi}_X\colon A\to \cT_X$ 
induces an isomorphism on the $K$-groups, 
it suffices to show that the $*$-ho\-mo\-mor\-phism 
$\bar{\pi}_X\colon A\to T$ 
induces an isomorphism on the $K$-groups. 
This can be shown by applying the discussion above 
to the non-degenerate \Cc in the next lemma. 

\begin{lemma}\label{non-deg}
Let us set $X'=\cspa (\varphi_X(A)X)$ which is a 
non-degenerate \Cc over $A$. 
Then there exists an isomorphism $\rho\colon \cT_{X'}\to T$ 
such that $\rho\circ \bar{\pi}_{X'}=\bar{\pi}_X$. 
\end{lemma}

\begin{proof}
Let us set $\pi=\bar{\pi}_X$ 
and define a linear map $t\colon X'\to \cT_X$ 
as the restriction of $\bar{t}_X$ to $X'$. 
It is easy to see that the pair $(\pi,t)$ 
is a representation of $X'$. 
Hence we have a $*$-ho\-mo\-mor\-phism 
$\rho\colon \cT_{X'}\to \cT_X$. 
It is clear that the gauge action of $\cT_X$ 
is a gauge action for the representation $(\pi,t)$. 
It is also clear that 
$\{a\in A\mid \pi(a)\in\psi_t(\cK(X'))\}=0$. 
Hence $\rho$ is injective by Theorem \ref{GIUTT}. 
Finally, it is not difficult to see that 
the image of $\rho$ is $T$. 
\end{proof}

This completes the proof of Proposition \ref{K(TX)}.

\end{document}